\documentclass[a4paper,11pt]{article}
\usepackage{bm,mathrsfs,amsmath,amsthm,amssymb,eepic,amscd,longtable,array,graphicx,mathdots}
\newcommand{\nc}{\newcommand}
\nc{\rnc}{\renewcommand}
\nc{\nn}{\nonumber}
\nc{\der}{{\partial}}
\rnc{\Im}{{\rm{Im}\,}}
\rnc{\Re}{{\rm{Re}\,}}
\nc{\db}{\displaybreak[0]\\}
\nc{\bra}{\langle}
\nc{\ket}{\rangle}
\nc{\bs}{\boldsymbol}

\DeclareMathOperator{\Tr}{Tr}

\DeclareMathOperator{\End}{End}

\newtheorem{theorem}{Theorem}[section]
\newtheorem{lemma}[theorem]{Lemma}

\newtheorem{proposition}[theorem]{Proposition}

\theoremstyle{definition}
\newtheorem{definition}[theorem]{Definition}

\numberwithin{equation}{section}

\numberwithin{equation}{section}

\begin{document}%
%
\title{Dual wavefunction of the Felderhof model}

\author{
Kohei Motegi \thanks{E-mail: kmoteg0@kaiyodai.ac.jp}
\\\\
{\it Faculty of Marine Technology, Tokyo University of Marine Science and Technology,}\\
 {\it Etchujima 2-1-6, Koto-Ku, Tokyo, 135-8533, Japan} \\
\\\\
\\
}

\date{\today}

\maketitle

\begin{abstract}
We study the Felderhof free-fermion six-vertex model,
whose wavefunction recently turned out to possess rich combinatorial
structure of the Schur polynomials.
We investigate the dual version of the wavefunction
in this paper, which seems to be a harder object to analyze.
We evaluate the dual wavefunction in two ways.
First, we give the exact correspondence between
the dual wavefunction and the Schur polynomials,
for which two proofs are given.
Next, we make a microscopic analysis and
express the dual wavefunction in terms of
strict Gelfand-Tsetlin pattern.
As a consequence of these two ways of evaluation
of the dual wavefunction, we obtain a dual version of the
Tokuyama combinatorial formula for the Schur polynomials.
We also give a generalization of the correspondence between the
dual wavefunction of the Felderhof model
and the factorial Schur polynomials.
\end{abstract}

\section{Introduction}
Integrable lattice models \cite{Bethe,FST,Baxter,KBI}
in mathematical physics
have played important roles in the developments
of algebras, combinatorics and representation theory.
One of the most fundamental models in integrable lattice models
is the six-vertex models \cite{Lieb,Reshetikhin}.
The most famous six-vertex model is the one
whose $L$-operator has the quantum group \cite{Dr,J} $U_q(sl_2)$ symmetry.
The corresponding one-dimensional integrable quantum spin chain
for this two-dimensional six-vertex model is the Heisenberg XXZ chain.

A less well-known six-vertex model is the
Felderhof free-fermion model \cite{Felderhof},
which can be regarded as a free fermion model
in an external field.
It was found some time before that the Felderhof model
has also quantum group symmetry
\cite{Mu,DA}. The corresponding representation has a property
that the quantum group parameter $q$ must
be roots of unity for the representation
to be finite-dimensional.
A special class of partition functions
called the domain wall boundary partition function
was also evaluated for the case of the Felderhof model
\cite{FCWZ} in the past.

However, it was only found in recent years that
the Felderhof model has rich mathematical structures
related with the combinatorial representation theory
of Schur polynomials.
One of the striking facts found \cite{BBF} was that
the Tokuyama formula \cite{To,OkTo}, which is a
one-parameter deformation of the 
Weyl character formula,
is naturally realized as wavefunction of
the Felderhof model.
The wavefunction is a special class of
partition function,
which can be called as an off-shell Bethe vector
since it becomes the Bethe eigenvectors of the
corresponding one-dimensional spin chain
when the Bethe ansatz equation is imposed
on the spectral parameters.
In this case, the wavefunction is sometimes called
as the on-shell Bethe vector.
However, we do not impose the Bethe ansatz equation
on the spectral parameters in this paper,
i.e., the parameters are free variables.

Besides the spectral parameter, one can introduce at least
one free parameter in the $L$-operator of the Felderhof model,
which turns out to play the role of the deformation parameter
in the Tokuyama formula for the Schur polynomials.
The parameter for the deformation can be
interpreted as a free paramater which can be
introduced when constructing a finite-dimensional
representation space of a quantum group
when $q$ is fixed at roots of unity.
Since the $L$-operator is constructed as
an intwertwiner of tensor product
of two representation spaces, one can in fat introduce at least
two free parameters,
one in the auxiliary space,
and another in the quantum space.
The parameters can
in principle be different for different auxiliary and quantum spaces. 
For the Tokuyama formula to be realized,
all the parameters are set to be equal in the auxiliary spaces, and
all are zero in the quantum spaces \cite{BBF}.
Keeping all the parameters in the quantum spaces non-zero and independent,
it was found that the wavefunction gives the factorial Schur polynomials
\cite{BMN}.
The Tokuyama formula for the Schur polynomials can be understood
as a consequence of the evaluation of the wavefunction
in two ways.
One by expressing it as a product of a one-parameter deformation
of the Vandermonde determinant and the Schur polynomials,
and another one by making a microscopic analysis and
deirve an expression using the strict Gelfand-Tsetlin pattern.
The Tokuyama formula is a consequence of the two evaluations
for the same object.
This understanding \cite{BBF} opened a new doorway
to the combinaotial representation theory of symmetric polynomials
via the Felderhof free-fermion model.

In this paper, we study the dual wavefunction
of the Felderhof model,
and study the combinatorics of the Schur polynomials
by analyzing the dual wavefunction,
a similar object but seems harder to analyze than the original wavefunction.
The dual wavefunction was
evaluated for the special case $t=1$ of the deformation parameter
\cite{BBF,BMN}, which was obtained by transforming the original wavefuncion
to the dual wavefunction by symmetry arguments.
We want the exact evaluation when the deformation parameter is generic,
since this free parameter plays the role of refining the combinatorics
of the Schur polynomials.
We evaluate the dual wavefunction in two ways
and obtain a combinatorial formula for the Schur polyomials.
First, we analyze the dual wavefunction directly,
and show the correspondence between the Schur polynomials.
We give two proofs for this correspondence,
one by using the arguments which is slightly more complicated than,
but the same with the one given in \cite{BBF}.
Another proof is a modern statistical mechanical approach,
which combines the matrix product method \cite{GMmat,KM}
and the Izergin-Korepin method of analysis
on the domain wall boundary partition function \cite{Ko,Iz}.
We next give a microscopic analysis of the dual wavefunction.
By calculating the matrix elements of a single $B$-operator,
we derive an expresstion of the dual wavefunction
in terms of the strict Gelfand-Tsetlin pattern.
By comparing the two evaluations of the dual wavefunction,
we derive a dual version of the Tokuyama-type formula
for the Schur polynomials.

This paper is organized as follows.
We introduce the Felderhof model in section 2 and
review the relation between the wavefunction and
the Schur polynomials in section 3. In sections 4 and 5,
we introduce the dual wavefunction,
and show the relation with the Schur polynomials
by giving two different proofs.
In section 5,
we evaluate the dual wavefunction based on the
calculation of the matrix elements of a single $B$-operator,
and express in terms of the strict Gelfand-Tsetlin pattern.
Combining the obtained expression with the one proved in sections 4 and 5,
we give a combinatorial formula which can be regarded
as a dual version of the Tokuyama formula.
We give a generalization of the correspondence between
the dual wavefunction of a generalization of the Felderhof model
and the factorial Schur polynomials in section 6.
Section 7 is devoted to the conclusion.

\section{Felderhof model}
We introduce the Felderhof model in this section,
and review
the results on the relation between the wavefunction and the Schur polynomials
in the next section.
We use the $L$-operator in \cite{BBF} which is
best suited for the study of the combinatorics of the Schur polynomials,
since the Tokuyama formula is exactly realized as the wavefunction
constructed from this $L$-operator.
More generic or gauge-transformed
ones can be found in \cite{Mu,DA,FCWZ} for example.
We also use the terminology of the quantum inverse scattering method
or the algebraic Bethe ansatz,
which is one of the most fundamental methods for the analysis
of quantum integrable models.

The most fundamental objects in integrable lattice models
are the $R$-matrix and the $L$-operator.
For the case of the Felderhof model,
the $R$-matrix is given by
\begin{eqnarray}
R_{ab}(z,t)=\left( 
\begin{array}{cccc}
1+tz & 0 & 0 & 0 \\
0 & t(1-z) & t+1 & 0 \\
0 & (t+1)z & z-1 & 0 \\
0 & 0 & 0 & z+t
\end{array}
\right), \label{rmatrix}
\end{eqnarray}
acting on the tensor product $W_a \otimes W_b$
of the complex two-dimensional space $W_a$.
Let us denote the orthonormal basis of $W_a$ and its dual as
$\{|0 \rangle_a, |1 \rangle_a \}$ and $\{{}_a \langle 0|, {}_a \langle 1|\}$,
and the matrix elements of the $R$-matrix as
$
{}_a \langle \gamma | {}_b \langle \delta | R_{a b}(z,t)
|\alpha \rangle_a | \beta \rangle_b=[R(z,t)]_{\alpha \beta}^{\gamma \delta}
$. The matrix elements of the $R$-matrix are explicitly given as
\begin{align}
{}_a \langle 0| {}_b \langle 0 | R_{a b}(z,t)
|0 \rangle_a | 0 \rangle_b&=1+tz, \\
{}_a \langle 0| {}_b \langle 1 | R_{a b}(z,t)
|0 \rangle_a | 1 \rangle_b&=t(1-z), \\
{}_a \langle 0| {}_b \langle 1 | R_{a b}(z,t)
|1 \rangle_a | 0 \rangle_b&=t+1, \\
{}_a \langle 1| {}_b \langle 0 | R_{a b}(z,t)
|0 \rangle_a |1 \rangle_b&=(t+1)z, \\
{}_a \langle 1| {}_b \langle 0 | R_{a b}(z,t)
|1 \rangle_a | 0 \rangle_b&=z-1, \\
{}_a \langle 1| {}_b \langle 1 | R_{a b}(z,t)
|1 \rangle_a | 1 \rangle_b&=z+t.
\end{align}

\begin{figure}[ht]
\includegraphics[width=15cm]{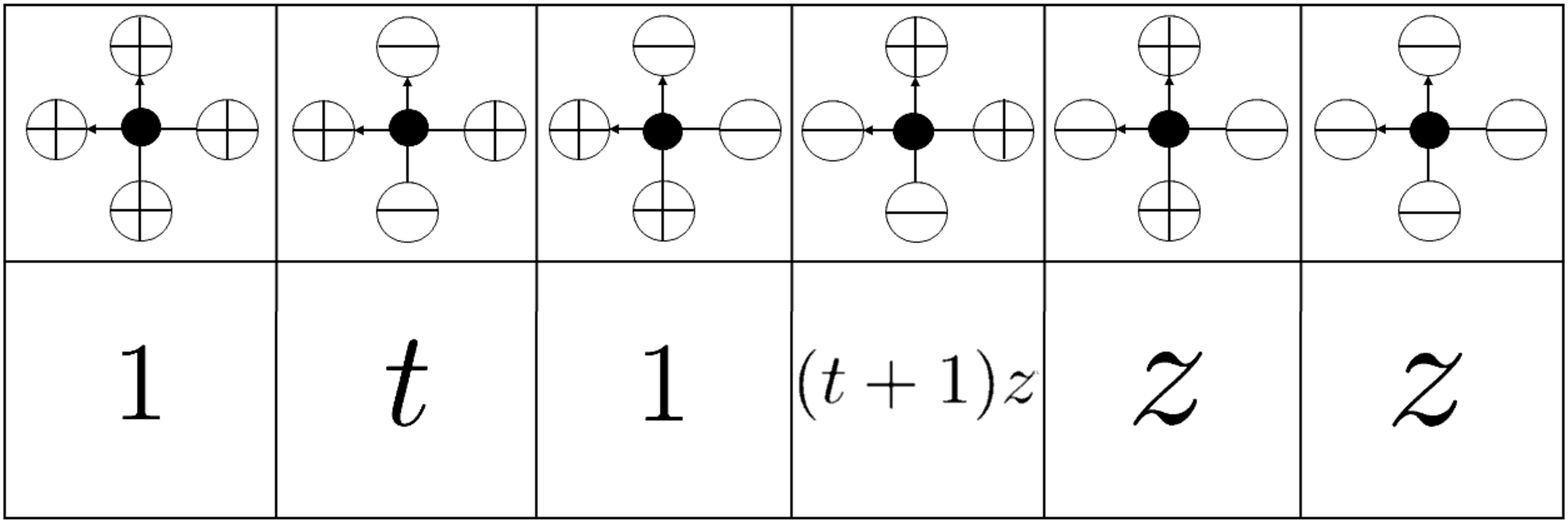}
\caption{The $L$-operator \eqref{loperator}. The (dual) state $| 0 \rangle$ ($\langle 0 |$) is represented as $\oplus$,
while the (dual) state $| 1 \rangle$ ($\langle 1 |$)
is represented as $\ominus$, following the pictorial description of \cite{BBF}.}
\label{pictureloperator}
\end{figure}

The $L$-operator of the Felderhof model is given by
\begin{eqnarray}
L_{aj}(z,t)=\left( 
\begin{array}{cccc}
1 & 0 & 0 & 0 \\
0 & t & 1 & 0 \\
0 & (t+1)z & z & 0 \\
0 & 0 & 0 & z
\end{array}
\right), \label{loperator}
\end{eqnarray}
acting on the tensor product $W_a \otimes \mathcal{F}_j$ of the space
$W_a$ and the two-dimensional Fock space at the
$j$th site $\mathcal{F}_j$.
We also denote the orthonormal basis of $\mathcal{F}_j$ and its dual as
$\{|0 \rangle_j, |1 \rangle_j \}$ and $\{{}_j \langle 0|, {}_j \langle 1|\}$,
and the matrix elements of the $L$-operator as
$
{}_a \langle \gamma| {}_j \langle \delta | L_{a j}(z,t)
|\alpha \rangle_a | \beta \rangle_j=[L(z,t)]_{\alpha \beta}^{\gamma \delta}
$.
The matrix elements of the $L$-operator are explicitly written as
(see Figure \ref{pictureloperator} for a pictorial description)
\begin{align}
{}_a \langle 0| {}_j \langle 0 | L_{a j}(z,t)
|0 \rangle_a | 0 \rangle_j&=1, \\
{}_a \langle 0| {}_j \langle 1 | L_{a j}(z,t)
|0 \rangle_a | 1 \rangle_j&=t, \\
{}_a \langle 0| {}_j \langle 1 | L_{a j}(z,t)
|1 \rangle_a | 0 \rangle_j&=1, \\
{}_a \langle 1| {}_j \langle 0 | L_{a j}(z,t)
|0 \rangle_a |1 \rangle_j&=(t+1)z, \\
{}_a \langle 1| {}_j \langle 0 | L_{a j}(z,t)
|1 \rangle_a | 0 \rangle_j&=z, \\
{}_a \langle 1| {}_j \langle 1 | L_{a j}(z,t)
|1 \rangle_a | 1 \rangle_j&=z.
\end{align}
The $R$-matrices and
the $L$-operators have origins in statistical physics,
and $| 0 \rangle$ or its dual $\langle 0|$
can be regarded as a hole state,
while $| 1 \rangle$ or its dual $\langle 1|$
can be interpretted as a particle state
from the point of view of statistical physics.
We use the terms hole states and particle states
to describe states constructed from
$| 0 \rangle$, $\langle 0|$, $| 1 \rangle$ and $\langle 1|$
from now on since they are convenient for the description
of the states.
We also remark that in the language of the quantum inverse scattering method,
the Fock spaces $W_a$ and $\mathcal{F}_j$
are usually called the auxiliary and quantum spaces, respectively.

The $R$-matrix \eqref{rmatrix}
and $L$-operator \eqref{loperator} satsify the Yang-Baxter relation
\begin{align}
R_{ab}(z_1/z_2,t)L_{aj}(z_1,t)L_{bj}(z_2,t)
=L_{bj}(z_2,t)L_{aj}(z_1,t)R_{ab}(z_1/z_2,t), \label{RLL}
\end{align}
acting on $W_a \otimes W_b \otimes V_j$.
We remark that this $RLL$ relation \eqref{RLL}
can be regarded as a special case of the
generalized Yang-Baxter relation for a more general
$R$-matrix \cite{Mu,DA,FCWZ}.
The $R$-matrix \eqref{rmatrix} and the $L$-operator
\eqref{loperator} in this section can be regarded as
different specializations of the general $R$-matrix
from this viewpoint.
One advantages of the point of view from the quantum group
used was that one can systematically generalize the Felderhof model
to higher-dimensional representations \cite{DA}.

\begin{figure}[h]
\includegraphics[width=15cm]{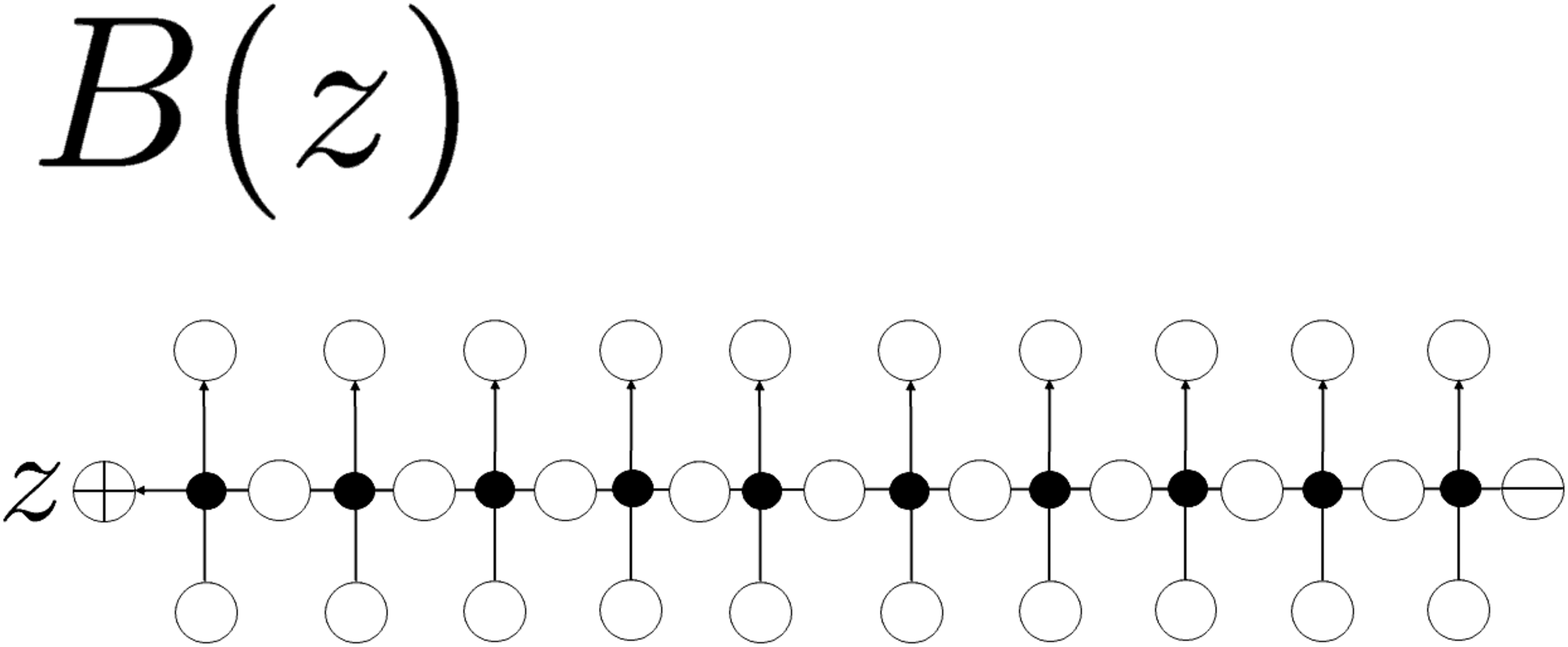}
\caption{The $B$-operator $B(z)$, which is a matrix element of
the monodromy matrix $T_a(z)$. The $B$-operator is
$2^M \times 2^M$ matrix-valued.
The leftmost state on the horizontal line
(auxiliary space) is fixed as $\oplus$ ($_{a} \langle 0 |$),
whereas the rightmost state is fixed as $\ominus$ ($| 1 \rangle_a$).}
\label{pictureBoperator}
\end{figure}

From the $L$-operator, we construct the monodromy matrix
\begin{align}
T_{a}(z)=L_{a M}(z,t) \cdots L_{a 1}(z,t)=
\begin{pmatrix}
A(z) & B(z)  \\
C(z) & D(z)
\end{pmatrix}_{a}
,
\label{monodromy1}
\end{align}
which acts on $W_a \otimes (\mathcal{F}_1\otimes\dots\otimes 
\mathcal{F}_{M})$.
The intertwining relation between the monodromy matrices
\begin{align}
R_{ab}(z_1/z_2,t)T_{a}(z_1)T_{b}(z_2)
=T_{b}(z_2)T_{a}(z_1)R_{ab}(z_1/z_2,t), \label{RTT}
\end{align}
follow from the $RLL$ relation \eqref{RLL}.
One of the elements of \eqref{RTT} is the commutation relations
between the $B$ operators
\begin{align}
(1+tz_1/z_2)B(z_1)B(z_2)=B(z_2)B(z_1)(z_1/z_2+t). \label{Boperatorscommutation}
\end{align}
Note that unlike the one constructed from the usual $U_q(sl_2)$ $R$-matrix,
the $B$-operators created from the Felderhof model \eqref{loperator}
do no simply commute, i.e., they produce extra factors.
See Figure \ref{pictureBoperator} for a graphical
description of the $B$-operator.

\section{Wavefunction and Schur polynomials}
We introduce the wavefunction which is a special class of partition functions,
and review how it is related with the
Schur polynomials defined below.
\begin{definition}
The Schur polynomials is defined to be the following determinant:
\begin{align}
s_\lambda(\{ z \}_N)=
   \frac{\mathrm{det}_N(z_j^{\lambda_k+N-k})}
        {\prod_{1 \le j < k \le N}(z_j-z_k)},
 \label{Schur}
\end{align}
where $\{ z \}_N=\{z_1,\dots,z_N \}$ is a set of variables
and $\lambda$ denotes a Young diagram
$\lambda=(\lambda_1,\lambda_2,\dots,\lambda_N)$
with weakly decreasing non-negative integers
$\lambda_1 \ge \lambda_2 \ge \cdots \ge \lambda_N \ge 0$.
\end{definition}
Before introducing the wavefunction, we first define
the arbitrary $N$-particle state \\
$|\Psi(z_1,\dots,z_N) \ket$
 with $N$ spectral parameters
$\{ z \}_N=\{ z_1,\dots,z_N \}$
as a state constructed by a multiple action
of $B$ operator on the vacuum state 
$|\Omega \ket:= | 0^{M} \ket:=|0\ket_1\
\otimes \dots \otimes |0\ket_{M}$
\begin{align}
|\Psi(z_1,\dots,z_N) \ket=B(z_1) \cdots B(z_N)| \Omega \ket.
\label{statevector}
\end{align}
\eqref{statevector} is usually called as the off-shell Bethe vector
(off-shell means that we do not assume the spectral parameters
satsify the Bethe ansatz equations).
Note that due to the commutation relation between the $B$-operators
\eqref{Boperatorscommutation},
the ordering of the $B$-operators in the off-shell Bethe vector
\eqref{statevector} is important for the Felderhof model.

Next, we introduce the wavefunction
$\bra x_1 \cdots x_N | \Psi(z_1,\dots,z_N) \ket$
as the overlap between an arbitrary off-shell 
$N$-particle state $|\Psi(z_1,\dots,z_N)\ket$ and 
the (normalized) state with an arbitrary particle configuration 
$|x_1 \cdots x_N\ket$ $(1 \le x_1<\dots<x_N \le M$), 
where $x_j$ denotes the positions of the particles. 
The particle configurations are explicitly defined as
\begin{align}
   \bra x_1 \cdots x_N|=\bra \Omega|\prod_{j=1}^N \sigma^+_{x_j},
\end{align}
where $\langle \Omega|:=\langle 0^{M}|:=
{}_1\bra 0|\otimes\dots \otimes{}_M\bra 0|$.
Here, we define $\sigma^+$ and $\sigma^-$ as
operators acting on the basis elements as
\begin{align}
&\sigma^+|1 \rangle=|0 \rangle, \ 
\sigma^+|0 \rangle=0, \ 
\langle 0|\sigma^+=\langle 1|, \
\langle 1|\sigma^+=0, 
\\
&\sigma^-|0 \rangle=|1 \rangle, \
\sigma^-|1 \rangle=0, \
\langle 1|\sigma^-=\langle 0|, \
\langle 0|\sigma^-=0.
\end{align}
The subscript $j$ of $\sigma_j^+$ or $\sigma_j^-$
indicates  that the operator acts on the space $\mathcal{F}_j$ as $\sigma^+$
or $\sigma^-$,
and as an idenitity on the other spaces.

Bump, Brubaker and Friedberg proved the following relation
between the wavefunction of the Felderhof model
and the Schur polynomials.
\begin{theorem} \label{TheoremBBF} \cite{BBF}
The wavefunction $\bra x_1 \cdots x_N | \Psi(z_1,\dots,z_N) \ket$
is expressed by the Schur polynomials as
\begin{align}
\bra x_1 \cdots x_N | \Psi(z_1,\dots,z_N) \ket
=\prod_{1 \le j < k \le N}(z_j+tz_k)s_\lambda(\{ z \}_N).
\end{align}
Here the Young diagram for the Schur polynomials correspond to
the particle configuration under the relation
$\lambda_j=x_{N-j+1}-N+j-1$, $j=1,\dots,N$.
\end{theorem}

The authors in \cite{BBF} moreover found that
the investigating the microscopic description of
the wavefunction of the Felderhof model
naturally leads to the Tokuyama formula
for the Schur polynomials,
which is a deformation of the Weyl character formula.
The idea is as follows.
First, we introduce
a strict Gelfand-Tsetlin pattern, which is a triangular array of integers
\begin{eqnarray}
\mathcal{T}=\left\{
\begin{array}{ccccccc}
a_{0,0} & & a_{0,1} & \cdots & a_{0,N-2} & & a_{0,N-1} \\
& a_{1,1} & & \cdots & & a_{1,N-1} & \\
& & \ddots & & \iddots & & \\
& & & a_{N-1,N-1} & & &
\end{array}
\right\},
\end{eqnarray}
in which the rows interlace $a_{i-1,j-1} \ge a_{i,j} \ge a_{i-1,j}$,
and the entries in horizontal rows are strictly decreasing.

For each strict Gelfand-Tsetlin pattern, we assign
the following weight
\begin{align}
G(\mathcal{T},\{ z \}_N)
=\prod_{i=1}^{N-1}
\prod_{j=i}^{N-1} \gamma(a_{i,j})
z_1^{d_0(\mathcal{T})-d_1(\mathcal{T})}
z_2^{d_1(\mathcal{T})-d_2(\mathcal{T})} \cdots
z_{N-1}^{d_{N-2}(\mathcal{T})-d_{N-1}(\mathcal{T})}
z_N^{d_{N-1}(\mathcal{T})}, \label{productofweights1}
\end{align}
where
$d_j(\mathcal{T})=\sum_{k=j}^{N-1}a_{j,k}$, $j=0,\dots,N-1$
is the sum of the entries of the strict Gelfand-Tsetlin pattern
in the $j$-th row, and $\gamma(a_{i,j})$ is defined as
\begin{align}
\gamma(a_{i,j})
=
\left\{
\begin{array}{ll}
t & a_{i,j}=a_{i-1,j-1}, \\
t+1 & a_{i-1,j} \neq a_{i,j} \neq a_{i-1,j-1}, \\
1 & \mathrm{otherwise}
\end{array}
\right. \label{gamma}
\end{align}
for pairs of integers $(i,j)$ satisfying
$1 \le i \le N-1, \ i \le j \le {N-1}$.

Investigating the inner states making nonzero contributions
to the wavefunction, one finds that the corresponding weight
for a fixed inner state,
which is the product of the matrix elements of the $L$-operators
of the inner states,
can be characterized by a strict Gelfand-Tsetlin pattern
with the top row fixed by the Young diagram as
$a_{0,j}=\lambda_{j+1}+N-j-1$.
The weight for each inner state is found to be given by
\eqref{productofweights1},
and the wavefunction can be expressed as a sum
of \eqref{productofweights1} for all strict Gelfand-Tsetlin patterns
with the top row fixed as $a_{0,j}=\lambda_{j+1}+N-j-1$.
Combining this microscopic analysis with
Theorem \ref{TheoremBBF}, one gets the following
combinatorial formula for the Schur polynomials.
\begin{theorem} \label{combinatorialformula} \cite{BBF}
We have the following combinatorial formula for the Schur polynomials
\begin{align}
&\prod_{1 \le j < k \le N}(z_j+tz_k)s_\lambda(\{ z \}_N) \nonumber \\
=&\sum_\mathcal{T} G(\mathcal{T},\{ z \}_N) \nonumber \\
=&\sum_\mathcal{T} \prod_{i=1}^{N-1}
\prod_{j=i}^{N-1} \gamma(a_{i,j})
z_1^{d_0(\mathcal{T})-d_1(\mathcal{T})}
z_2^{d_1(\mathcal{T})-d_2(\mathcal{T})} \cdots
z_{N-1}^{d_{N-2}(\mathcal{T})-d_{N-1}(\mathcal{T})}
z_N^{d_{N-1}(\mathcal{T})}, \label{combinatorialformulaequation}
\end{align}
where the sum is over all strict Gelfand-Tsetlin patterns
with the top row of the strict Gelfand-Tsetlin pattern
is fixed by the Young diagram as
$a_{0,j}=\lambda_{j+1}+N-j-1$.
\end{theorem}

\section{Dual wavefunction}
We now introduce the dual wavefunction, and
study the exact relation between it
and the Schur polynomials.
In this section,
we use the argument which is more slightly complicated than, but basically
the same with the one given in \cite{BBF}.
We analyze by another method based on a modern statistical physical
and quantum integrable techniques, which will be
given in the next section.

Before defining the dual wavefunction,
we introduce another type of arbitrary dual $N$-hole state
$\langle \Phi(z_1,\dots,z_N)|$ by a multiple action of $B$
operator on the dual particle occupied state
$\langle 1 \cdots M|:=\langle 1^{M}|:=
{}_1\bra 1|\otimes\dots \otimes{}_M\bra 1|$
\begin{align}
\langle \Phi(z_1,\dots,z_N)|
=\langle 1 \cdots M|B(z_1) \cdots B(z_N).
\end{align}
It is convenient to introduce a notation for
the state with an arbitrary hole configuration 
$|\overline{x_1} \cdots \overline{x_N} \ket$
$(1 \le \overline{x_1}<\dots< \overline{x_N} \le M$), 
where $\overline{x_j}$ denotes the positions of holes. 
Explicitly,
\begin{align}
|\overline{x_1} \cdots \overline{x_N} \ket
=\prod_{j=1}^N \sigma^+_{x_j}
(|1 \rangle_1 \otimes \cdots \otimes |1 \rangle_M).
\end{align}

\begin{figure}[ht]
\includegraphics[width=15cm]{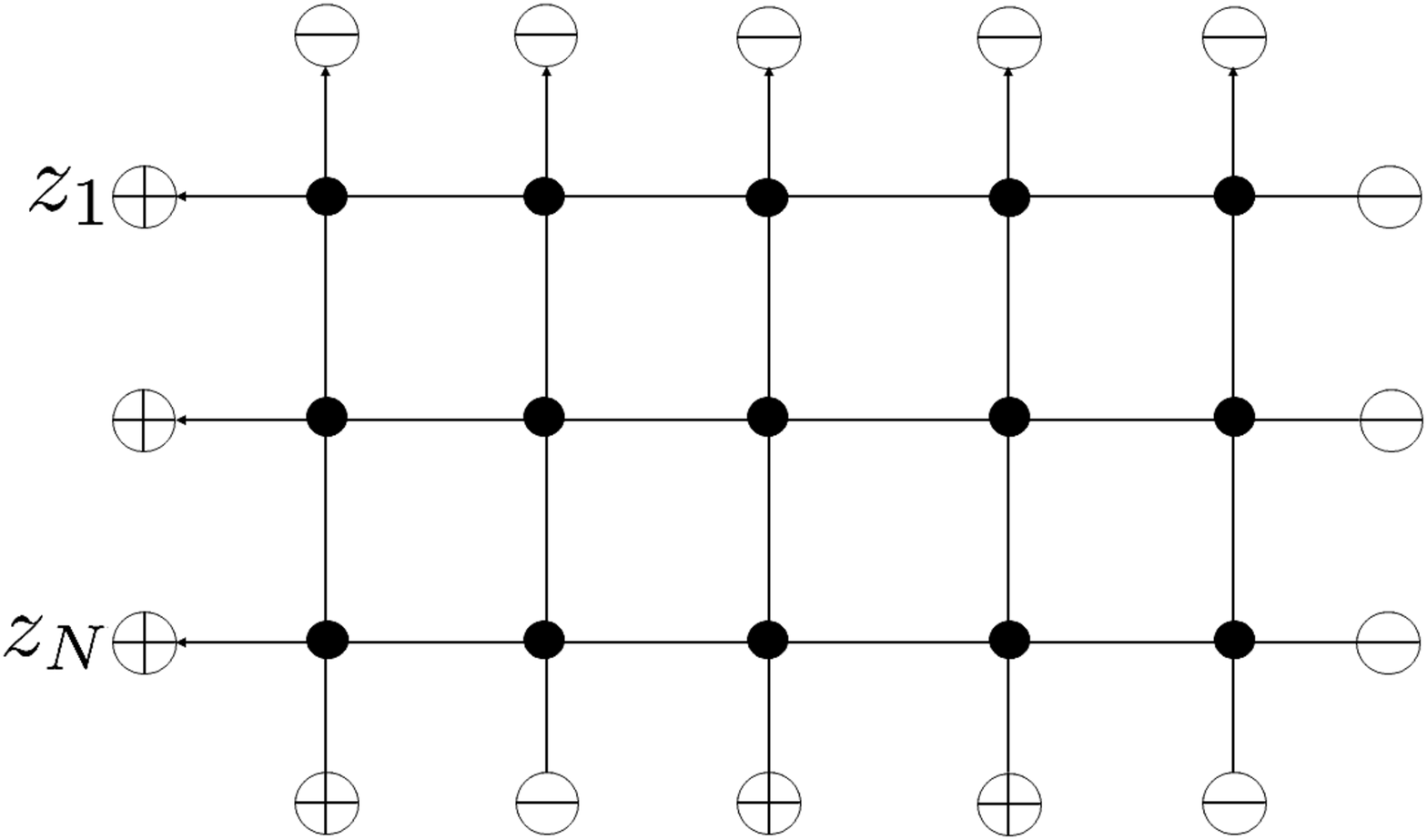}
\caption{The dual wavefunction
$\langle \Phi(z_1,\dots,z_N)|\overline{x_1} \cdots \overline{x_N} \ket$
for the case $M=5$, $N=3$, $(\overline{x_1},\overline{x_2},\overline{x_3})
=(2,3,5)$.}
\label{picturedualwavefunction}
\end{figure}

The dual wavefunction
$\langle \Phi(z_1,\dots,z_N)|\overline{x_1} \cdots \overline{x_N} \ket$
is defined as the overlap between the
arbitrary dual $N$-hole state
$\langle \Phi(z_1,\dots,z_N)|$
and hole configurations
$|\overline{x_1} \cdots \overline{x_N} \ket$
(see Figure \ref{picturedualwavefunction} for an example
of a graphical description of the dual wavefunction).

We show the following relation between the dual wavefunction
and the Schur polynomials.
\begin{theorem} \label{theoremdualandschur}
The dual wavefunction
$\langle \Phi(z_1,\dots,z_N)| \overline{x_1} \cdots \overline{x_N}
\rangle$ can be expressed by the Schur polynomials as
\begin{align}
\langle \Phi(z_1,\dots,z_N)| \overline{x_1} \cdots \overline{x_N}
\rangle=t^{N(M-N)} \prod_{1 \le j < k \le N}(z_j+t z_k)
s_{\overline{\lambda}} \Bigg( \Bigg\{ \frac{z}{t} \Bigg\}_N \Bigg).
\label{dualwavefunctionsandschur}
\end{align}
Here the Young diagram for the Schur polynomials correspond to
the particle configuration under the relation
$\overline{\lambda_j}=\overline{x_{N-j+1}}-N+j-1$, $j=1,\dots,N$,
and the symmetric variables are
$\displaystyle \Bigg\{ \frac{z}{t} \Bigg\}_N=
\Bigg\{ \frac{z_1}{t},\dots,\frac{z_N}{t} \Bigg\}$.
\end{theorem}
Before proving Theorem \ref{theoremdualandschur},
let us make some comments.
There is a factor $t^{N(M-N)}$ which depends on the number
of sites $M$ and the number of particles $N$ in the right hand side
of \eqref{dualwavefunctionsandschur}.
What makes things more complicated is that the
symmetric variables of the Schur polynomials
are $\displaystyle \Bigg\{ \frac{z}{t} \Bigg\}_N$,
which are the reasons why the relation
for the dual wavefunction
\eqref{dualwavefunctionsandschur} seems hard to find.
We actually first found this Theorem by using the statistical physical
method given in the next section.
One advantages of the proof given in this section following \cite{BBF}
is that the proof naturally lifts to the correspondence between
a generalization of the Felderhof model and the factorial Schur polynoimals.
\begin{proof}
We rewrite \eqref{dualwavefunctionsandschur}
by rescaling each $z_j$ to $t z_j$ \eqref{dualwavefunctionsandschur}
in the following form
\begin{align}
t^N
\langle 1 \cdots M|\Bigg(\frac{B(tz_1)}{t^M}\Bigg) \cdots
\Bigg(\frac{B(tz_N)}{t^M} \Bigg)| \overline{x_1} \cdots \overline{x_N}
\rangle=\prod_{1 \le j < k \le N} 
\Bigg(\frac{1}{t}z_j+z_k \Bigg)
s_{\overline{\lambda}} ( \{ z \}_N ). \label{equivalentone}
\end{align}
For giving a proof, it is convenient to
introduce the rescaled $L$-operator
\begin{eqnarray}
\widetilde{L}(z,t)=
\frac{1}{t}L(tz,t)=
\left( 
\begin{array}{cccc}
1/t & 0 & 0 & 0 \\
0 & 1 & 1/t & 0 \\
0 & (t+1)z & z & 0 \\
0 & 0 & 0 & z
\end{array}
\right), \label{loperator2}
\end{eqnarray}
and the rescaled monodromy matrix
\begin{align}
\widetilde{T}_{a}(z)=\widetilde{L}_{a M}(z,t) \cdots \widetilde{L}_{a 1}(z,t)=
\begin{pmatrix}
\widetilde{A}(z) & \widetilde{B}(z)  \\
\widetilde{C}(z) & \widetilde{D}(z)
\end{pmatrix}_{a}
.
\label{monodromy2}
\end{align}
Using these rescaled objects,
\eqref{equivalentone} can be expressed as
\begin{align}
t^N
\langle 1 \cdots M|\widetilde{B}(z_1) \cdots
\widetilde{B}(z_N)| \overline{x_1} \cdots \overline{x_N}
\rangle=\prod_{1 \le j < k \le N} 
\Bigg(\frac{1}{t}z_j+z_k \Bigg)
s_{\overline{\lambda}} (\{ z \}_N). \label{equivalenttwo}
\end{align}
Instead of proving \eqref{dualwavefunctionsandschur},
we show \eqref{equivalenttwo} since this is equivalent to
\eqref{dualwavefunctionsandschur} and is the expression
which one can use the argument given in \cite{BBF}.

We first show the following lemma.

\begin{lemma} \label{lemmaone}
\begin{align}
\prod_{1 \le j < k \le N} 
\Bigg(\frac{1}{t}z_j+z_k \Bigg)^{-1}
t^N
\langle 1 \cdots M|\widetilde{B}(z_1) \cdots
\widetilde{B}(z_N)| \overline{x_1} \cdots \overline{x_N}
\rangle,
\end{align}
does not depend on $t$.
\end{lemma}
\begin{proof}
We prove this lemma by showing the following properties
for \\
$
t^N
\langle 1 \cdots M|\widetilde{B}(z_1) \cdots
\widetilde{B}(z_N)| \overline{x_1} \cdots \overline{x_N}
\rangle
$: 
\begin{enumerate}
\item 
$
t^N
\langle 1 \cdots M|\widetilde{B}(z_1) \cdots
\widetilde{B}(z_N)| \overline{x_1} \cdots \overline{x_N}
\rangle
$
is a polynomial of $t^\prime:=t^{-1}$ with highest degree $N(N-1)/2$. 
\item
$
t^N
\langle 1 \cdots M|\widetilde{B}(z_1) \cdots
\widetilde{B}(z_N)| \overline{x_1} \cdots \overline{x_N}
\rangle
$
has $t^\prime z_j+z_k$, $1 \le j < k \le N$ as factors. 
\end{enumerate}

We first show
$
\mathrm{deg}_{t^\prime}
(t^N
\langle 1 \cdots M|\widetilde{B}(z_1) \cdots
\widetilde{B}(z_N)| \overline{x_1} \cdots \overline{x_N}
\rangle)
\le N(N-1)/2
$ by induction on $N$.
The case $N=1$ follows as an special case of the general fact \\
$
0 \le \mathrm{deg}_{t^\prime}
(t \langle \overline{x_1} \cdots \overline{x_N} |\widetilde{B}(z)| \overline{y_1}
\cdots \overline{y_{N+1}} \rangle) \le N
$
which can be seen easily from the definition of
the rescaled $L$-operator $\widetilde{L}(z,t)$.
Next, let us assume $
\mathrm{deg}_{t^\prime}
(t^N
\langle 1 \cdots M|\widetilde{B}(z_1) \cdots
\widetilde{B}(z_N)| \overline{x_1} \cdots \overline{x_N}
\rangle)
\le N(N-1)/2
$.
One can see
$
\mathrm{deg}_{t^\prime}
(t^{N+1}
\langle 1 \cdots M|\widetilde{B}(z_1) \cdots
\widetilde{B}(z_{N+1})| \overline{y_1} \cdots \overline{y_{N+1}}
\rangle)
\le (N+1)N/2
$
by  combining the assumption
$
\mathrm{deg}_{t^\prime}
(t^N
\langle 1 \cdots M|\widetilde{B}(z_1) \cdots
\widetilde{B}(z_N)| \overline{x_1} \cdots \overline{x_N}
\rangle)
\le N(N-1)/2
$, \\
the fact
$
0 \le \mathrm{deg}_{t^\prime}
(t \langle \overline{x_1} \cdots \overline{x_N} |\widetilde{B}(z)| \overline{y_1}
\cdots \overline{y_{N+1}} \rangle) \le N
$
and the decomposition
\begin{align}
&t^{N+1}
\langle 1 \cdots M|\widetilde{B}(z_1) \cdots
\widetilde{B}(z_{N+1})| \overline{y_1} \cdots \overline{y_{N+1}}
\rangle \nonumber \\
=&\sum_{\{ \overline{x} \}}
(t^N
\langle 1 \cdots M|\widetilde{B}(z_1) \cdots
\widetilde{B}(z_{N})| \overline{x_1} \cdots \overline{x_{N}} \rangle)
(t \langle \overline{x_1} \cdots \overline{x_N} |\widetilde{B}(z_{N+1})| \overline{y_1}
\cdots \overline{y_{N+1}} \rangle).
\end{align}

Next we show Property 2.
The commutation relation
\eqref{Boperatorscommutation}
can be rewritten as the following commutation relation
between the rescaled $B$-operators
\begin{align}
(z_1+t^\prime z_2)\widetilde{B}(z_1)\widetilde{B}(z_2)
=\widetilde{B}(z_2)\widetilde{B}(z_1)(t^\prime z_1+z_2).
\label{rescaledcommutation}
\end{align}
Applying the commutation relation
\eqref{rescaledcommutation} repeatedly,
one gets the following equality
\begin{align}
&t^N
\langle 1 \cdots M|\widetilde{B}(z_1) \cdots
\widetilde{B}(z_N)| \overline{x_1} \cdots \overline{x_N}
\rangle \prod_{1 \le j < k \le N} 
(z_j+t^\prime z_k) \nonumber \\
=&t^N
\langle 1 \cdots M|\widetilde{B}(z_N) \cdots
\widetilde{B}(z_1)| \overline{x_1} \cdots \overline{x_N}
\rangle \prod_{1 \le j < k \le N} 
(t^\prime z_j+z_k)
. \label{discussionequality}
\end{align}
Note that in the equality
\eqref{discussionequality}, the factors
$t^N
\langle 1 \cdots M|\widetilde{B}(z_1) \cdots
\widetilde{B}(z_N)| \overline{x_1} \cdots \overline{x_N}
\rangle$, \\
$\prod_{1 \le j < k \le N} 
(z_j+t^\prime z_k)$,
$t^N
\langle 1 \cdots M|\widetilde{B}(z_N) \cdots
\widetilde{B}(z_1)| \overline{x_1} \cdots \overline{x_N}
\rangle$ and $\prod_{1 \le j < k \le N} 
(t^\prime z_j+z_k)
$ are polynomials of $t^\prime$.
From this fact and that $\prod_{1 \le j < k \le N} 
(z_j+t^\prime z_k)$ is not divided by $\prod_{1 \le j < k \le N} 
(t^\prime z_j+z_k)$,
one can see $t^N
\langle 1 \cdots M|\widetilde{B}(z_1) \cdots
\widetilde{B}(z_N)| \overline{x_1} \cdots \overline{x_N}
\rangle$ is divided by $\prod_{1 \le j < k \le N} 
(t^\prime z_j+z_k)$.

From Property 2, we have
$
\mathrm{deg}_{t^\prime}
(t^N
\langle 1 \cdots M|\widetilde{B}(z_1) \cdots
\widetilde{B}(z_N)| \overline{x_1} \cdots \overline{x_N}
\rangle)
\ge N(N-1)/2
$.
Together with
$
\mathrm{deg}_{t^\prime}
(t^N
\langle 1 \cdots M|\widetilde{B}(z_1) \cdots
\widetilde{B}(z_N)| \overline{x_1} \cdots \overline{x_N}
\rangle)
\le N(N-1)/2
$ which is proved before, we have Property 1.

\end{proof}

From Lemma \ref{lemmaone}, one sees that to study the wavefunction,
it is enough to examine a particular value of $t$.
The case when $t=-1$ in which the six-vertex model
reduces to a five-vertex model 
\begin{eqnarray}
\widetilde{L}(z,-1)=
-L(-z,-1)=
\left( 
\begin{array}{cccc}
-1 & 0 & 0 & 0 \\
0 & 1 & -1 & 0 \\
0 & 0 & z & 0 \\
0 & 0 & 0 & z
\end{array}
\right),
\end{eqnarray}
is easy to examine,
and we have the following relation.

\begin{lemma}
We have
\begin{align}
\prod_{1 \le j < k \le N} 
\Bigg(\frac{1}{t}z_j+z_k \Bigg)^{-1}
t^N
\langle 1 \cdots M|\widetilde{B}(z_1) \cdots
\widetilde{B}(z_N)| \overline{x_1} \cdots \overline{x_N}
\rangle \Bigg|_{t=-1}
=s_{\overline{\lambda}}(\{ z \}_N). \label{forproof}
\end{align}
\end{lemma}
\begin{proof}
To prove the Lemma is equivalent to show
\begin{align}
\langle 1 \cdots M|\widetilde{B}(z_1) \cdots
\widetilde{B}(z_N)| \overline{x_1} \cdots \overline{x_N}
\rangle|_{t=-1}=(-1)^N
\prod_{1 \le j < k \le N}
(-z_j+z_k) s_{\overline{\lambda}} (\{ z \}_N).
\end{align}

To show this, we first note that the matrix elements
of a single $B$-operator is given by
\begin{align}
\langle \overline{x_1} \cdots \overline{x_{k-1}}
|\widetilde{B}(z)| \overline{y_1}
\cdots \overline{y_{k}} \rangle
=&(-1)^k (-1)^{j-1} z^{\overline{y_j}-1},
\label{matrixelementst=-1}
\end{align}
when the hole configurations $\{ \overline{x} \}$ and $\{ \overline{y} \}$
satisfy \\
$\overline{x_1}=\overline{y_1}, \cdots,
\overline{x_{j-1}}=\overline{y_{j-1}}$,
$\overline{x_{j}}=\overline{y_{j+1}}, \cdots,
\overline{x_{k-1}}=\overline{y_{k}}$ for some $j$,
and 0 otherwise.

Since the matrix elements of a single $B$-operator are essentially the same
with the ones for the original wavefunction at $t=-1$ in \cite{BBF}
except the sign $(-1)^k$ 
(we also have to translate the hole configurations to
particle configurations), the same argument can be appplied.
One observes that
the number of the inner states whose weights gives non-zero contributions
to the dual wavefunction is $N!$.
The weight of each nonvanishing inner state
corresponds to one term $(-1)^\sigma
\prod_{j=1}^N z_j^{\overline{\lambda_{\sigma}(j)}+N-\sigma(j)}$
of the determinant expansion
of the numerator \eqref{Schur}
of the Schur polynomials
multipled by the extra factor $(-1)^{N(N+1)/2}$.
The factor $(-1)^{N(N+1)/2}$
appears since the dual wavefunction is constructed from
$N$ layers of $B$-operators, and
the $k$-th layer of the $B$-operator has
the extra sign $(-1)^k$ in the right hand side of
\eqref{matrixelementst=-1}, hence the total contribution
of $N$ layers of the $B$ operators gives the
extra sign $\prod_{k=1}^N (-1)^k=(-1)^{N(N+1)/2}$.
We further rewrite the extra factor $(-1)^{N(N+1)/2}$ as
\begin{align}
(-1)^{N(N+1)/2}=(-1)^N (-1)^{N(N-1)/2}
=(-1)^N
\frac{\prod_{1 \le j < k \le N}
(-z_j+z_k)}
{\prod_{1 \le j < k \le N}
(z_j-z_k)
},
\end{align}
to get
\begin{align}
&\langle 1 \cdots M|\widetilde{B}(z_1) \cdots
\widetilde{B}(z_N)| \overline{x_1} \cdots \overline{x_N}
\rangle|_{t=-1} \nonumber \\
=&
(-1)^{N(N+1)/2} \sum_{\sigma \in S_N} (-1)^\sigma
\prod_{j=1}^N z_j^{\overline{\lambda_{\sigma}(j)}+N-\sigma(j)}
\nonumber \\
=&(-1)^N
\frac{\prod_{1 \le j < k \le N}
(-z_j+z_k)}
{\prod_{1 \le j < k \le N}
(z_j-z_k)
}
\sum_{\sigma \in S_N} (-1)^\sigma
\prod_{j=1}^N z_j^{\overline{\lambda_{\sigma}(j)}+N-\sigma(j)}
\nonumber \\
=&(-1)^N
\prod_{1 \le j < k \le N}
(-z_j+z_k) s_{\overline{\lambda}} (\{ z \}_N)
.
\end{align}
\end{proof}
From Lemma \ref{lemmaone} and \eqref{forproof},
we have
\begin{align}
&\prod_{1 \le j < k \le N} 
\Bigg(\frac{1}{t}z_j+z_k \Bigg)^{-1}
t^N
\langle 1 \cdots M|\widetilde{B}(z_1) \cdots
\widetilde{B}(z_N)| \overline{x_1} \cdots \overline{x_N}
\rangle \nonumber \\
=&
\prod_{1 \le j < k \le N} 
\Bigg(\frac{1}{t}z_j+z_k \Bigg)^{-1}
t^N
\langle 1 \cdots M|\widetilde{B}(z_1) \cdots
\widetilde{B}(z_N)| \overline{x_1} \cdots \overline{x_N}
\rangle \Bigg|_{t=-1} \nonumber \\
=&s_{\overline{\lambda}}(\{ z \}_N),
\end{align}
which is exactly
\eqref{equivalenttwo}, hence
Theorem \ref{theoremdualandschur} is proved.
\end{proof}

\begin{figure}[h]
\includegraphics[width=11cm]{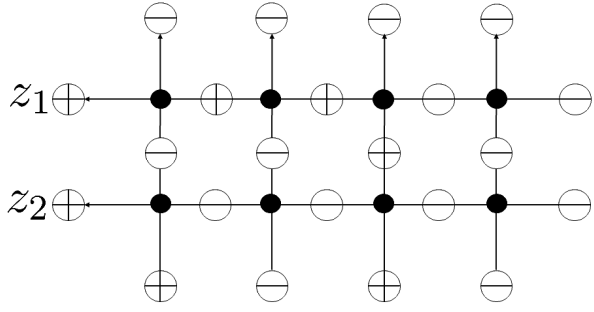}
\caption{An inner state of the dual wavefunction
$\langle \Phi(z_1,z_2)| \overline{x_1}=2, \overline{x_2}=4\rangle$
giving the weight $t^2 z_1 z_2^3$.}
\label{picturedualwavefunctioncontribution1}
\end{figure}

\begin{figure}[h]
\includegraphics[width=11cm]{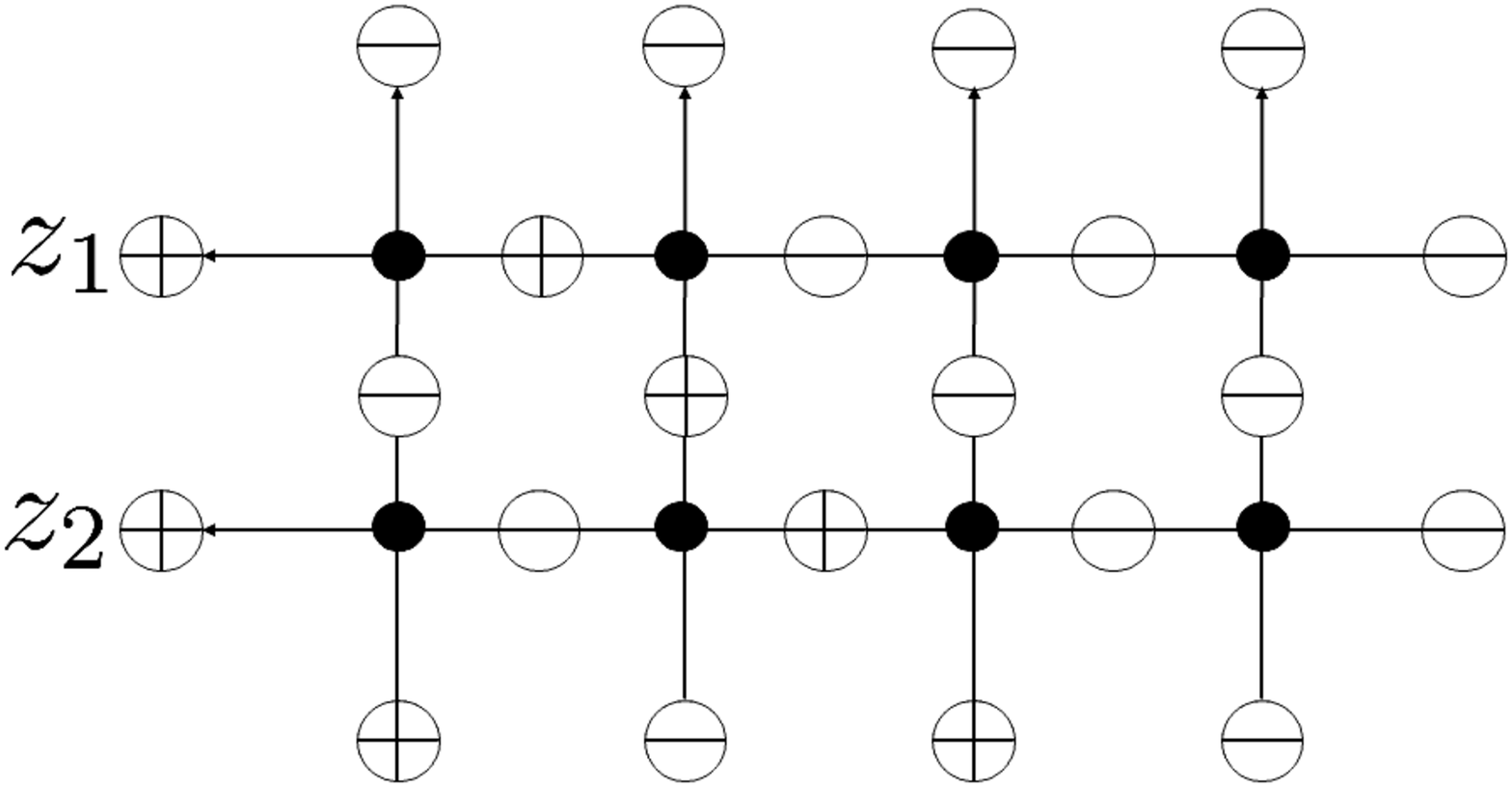}
\caption{An inner state of the dual wavefunction
$\langle \Phi(z_1,z_2)| \overline{x_1}=2, \overline{x_2}=4\rangle$
giving the weight $t(t+1) z_1^2 z_2^2$.}
\label{picturedualwavefunctioncontribution2}
\end{figure}

\begin{figure}[h]
\includegraphics[width=11cm]{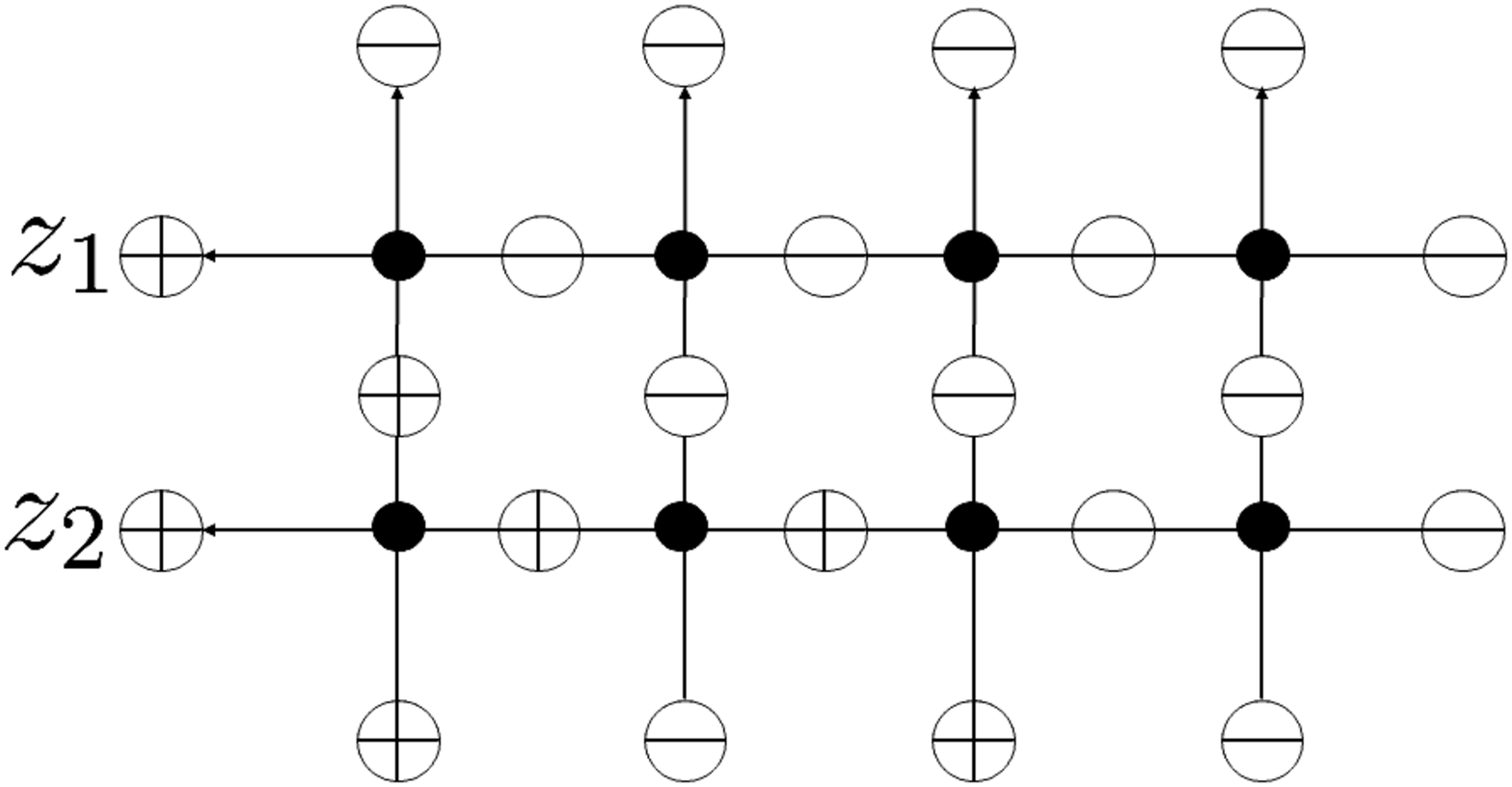}
\caption{An inner state of the dual wavefunction
$\langle \Phi(z_1,z_2)| \overline{x_1}=2, \overline{x_2}=4\rangle$
giving the weight $t z_1^3 z_2$.}
\label{picturedualwavefunctioncontribution3}
\end{figure}

{\bf Example}
Let us check Theorem \ref{theoremdualandschur} by an example
$M=4$, $N=2$ , $(\overline{x_1},\overline{x_2})=(2,4)$.
The corresponding Young diagram is
$\overline{\lambda}
=(\overline{\lambda}_1,\overline{\lambda_2})
=(\overline{x_2}-2,\overline{x_1}-1)=(4-2,2-1)=(2,1)$.
The left hand side of
\eqref{dualwavefunctionsandschur} can be calculated graphically
by noting that there are three inner states whose weights
make nonzero contributions
to the dual wavefunction
$\langle \Phi(z_1,z_2)| \overline{x_1}=2, \overline{x_2}=4\rangle$, 
which are given in Figures
\ref{picturedualwavefunctioncontribution1},
\ref{picturedualwavefunctioncontribution2}
and
\ref{picturedualwavefunctioncontribution3}.
From its graphical description
and using the data of the matrix elements of the
$L$-operator \eqref{loperator} and multiplying them,
one sees that
each of the configurations have weights
$t^2 z_1 z_2^3$, $t(t+1) z_1^2 z_2^2$
and $t z_1^3 z_2$.
Summing up the weights and noting that one can extract
$(z_1+t z_2)$ as an overall factor, we have
\begin{align}
t^2 z_1 z_2^3+t(t+1) z_1^2 z_2^2+t z_1^3 z_2
&=t^4 (z_1+t z_2) \Bigg( \Bigg(\frac{z_1}{t} \Bigg)^2
\Bigg(\frac{z_2}{t} \Bigg)
+\Bigg(\frac{z_1}{t} \Bigg)
\Bigg(\frac{z_2}{t} \Bigg)^2  
\Bigg) \nonumber \\
&=t^4 (z_1+t z_2)s_{(2,1)}\Bigg(\frac{z_1}{t}, \frac{z_2}{t} \Bigg),
\end{align}
which is nothing but the right hand side of
\eqref{dualwavefunctionsandschur}.

\section{Another proof}

We give another proof of Theorem \ref{theoremdualandschur}
by using a modern statistical mechanical method and
an analysis on a fundamental object in quantum integrable models,
i.e., we use the matrix product method
and the domain wall boundary partition function, as was
done in the case of the Grothendieck polynomials in \cite{MS2}
(see also \cite{MSW1} in which we demonstrate a proof of
Theorem \ref{TheoremBBF}
by using the same arguments given in this section).
We prove Theorem \ref{theoremdualandschur} as follows.
We first rewrite the dual wavefunction
into a matrix product form, following \cite{GMmat,KM}.
The matrix product form can be expressed as a determinant with some overall
factor which remains to be calculated. 
The information of the hole configuration 
$\{\overline{x_1},\overline{x_2},\dots,\overline{x_N} \}$
is encoded in the determinant.
On the other hand, the overall factor is independent of the
hole positions, and this factor can be determined by
considering the specific configuration: we 
explicitly evaluate the overlap of the consecutive hole configuration
(i.e. $\overline{x_j}=j$) whose evaluation
essentially reduces to that of the domain wall boundary partition function.

\begin{proof}
Let us begin to compute the wavefunction
$\bra \Phi(\{ z \}_N)| \overline{x_1} \cdots \overline{x_N} \ket$.
We first rewrite it into the matrix product representation.
With the help of graphical description,
one finds that the wavefunction can be written as
\begin{align}
\bra \Phi(\{ z \}_N)| \overline{x_1} \cdots \overline{x_N} \ket
=\Tr_{W^{\otimes N}}
\left[ 
\bra 1 \cdots M| \prod_{a=1}^N T_{a}(z_a) |
\overline{x_1} \cdots \overline{x_N} \ket P
\right],
\label{overlap}
\end{align}
where $P=| 1^N \rangle \langle 0^N |$
is an operator acting on the tensor product of auxiliary spaces
$W_1\otimes  \dots \otimes W_N$.
The trace here is also over the auxiliary spaces.

Changing the viewpoint of the products of the monodromy matrices, we have
\begin{align}
\prod_{a=1}^N T_{a}(z_a)
=\prod_{j=1}^M \mathcal{T}_j(\{z\}_N),
\end{align}
where $\mathcal{T}_j(\{z\}_N):= 
\prod_{a=1}^N L_{a j}(z_{a}) \in \End( W^{\otimes N} \otimes V_j)$
can be regarded as a monodromy matrix consisting of
$L$-operators acting on the same quantum space $V_j$
(but acting on different auxiliary spaces).  The monodromy matrix
is decomposed as
\begin{align}
\mathcal{T}_j(\{z\}_N)
&:=\begin{pmatrix}
\mathcal{A}_N (\{z\}_N) & \mathcal{B}_N(\{z\}_N) \\
\mathcal{C}_N(\{z\}_N) &  \mathcal{D}_N(\{z\}_N)
\end{pmatrix}_j ,
\label{decomp}
\end{align}
where the elements
($\mathcal{A}_N$, etc.) act on 
$W_1\otimes \dots \otimes W_N$.

The wavefunction \eqref{overlap} can then be rewritten by 
$\mathcal{T}_j(\{z\}_N)$ as
\begin{align}
\bra \Phi(\{ z \}_N)| \overline{x_1} \cdots \overline{x_N} \ket
&=\Tr_{W^{\otimes N}}
\left[ 
\bra 1 \cdots M| \prod_{j=1}^M \mathcal{T}_j(\{z\}_N) |
\overline{x_1} \cdots \overline{x_N} \ket P
\right] \nn \\
&=\Tr_{W^{\otimes N}}\left[ 
\mathcal{D}_N^{M-\overline{x_N}}
\mathcal{C}_N
\mathcal{D}_N^{\overline{x_N}-\overline{x_{N-1}}-1}
\dots\mathcal{C}_N\mathcal{D}_N^{\overline{x_2}-\overline{x_1}-1}
\mathcal{C}_N\mathcal{D}_N^{\overline{x_1}-1}
P
\right].
\label{reov}
\end{align}

For these operators, one finds the following recursive 
relations: 
\begin{align}
&\mathcal{D}_{n+1}(\{z\}_{n+1})
=\mathcal{D}_n(\{z\}_n) \otimes
\begin{pmatrix}
t & 0  \\
0 & z_{n+1}
\end{pmatrix} 
+\mathcal{C}_n(\{z\}_n) \otimes
\begin{pmatrix} 
0 & 0  \\
(1+t)z_{n+1} & 0
\end{pmatrix},
\label{reop1} \\
&\mathcal{C}_{n+1}(\{z\}_{n+1})
=\mathcal{D}_n(\{z\}_n) \otimes
\begin{pmatrix}
0 & 1  \\
0 & 0
\end{pmatrix}
+\mathcal{C}_n(\{z\}_n) \otimes
\begin{pmatrix}
1 & 0  \\
0 & z_{n+1}
\end{pmatrix}
,
\label{reop2}  
\end{align}
with the initial condition
\begin{align}
\mathcal{D}_1=
\begin{pmatrix}
t & 0 \\
0    &  z_1
\end{pmatrix}, \quad
\mathcal{C}_1=
\begin{pmatrix}
0 & 1 \\
0 & 0 
\end{pmatrix}.
\label{initialCD}
\end{align}

By using the recursive relations \eqref{reop1} and \eqref{reop2},
one sees that these operators satisfy the following simple algebra.
\begin{lemma}\label{algebra}
There exists a decomposition of $\mathcal{C}_n$ :
$\mathcal{C}_n=\sum_{j=1}^n \mathcal{C}_n^{(j)}$ such that
the following algebraic relations hold for $\mathcal{D}_n$ and $\mathcal{C}_n^{(j)}$:
\begin{align}
&\mathcal{C}_n^{(j)}\mathcal{D}_n
=\frac{z_j}{t}\mathcal{D}_n \mathcal{C}_n^{(j)}, \label{rel2} \\
&(\mathcal{C}_n^{(j)})^2=0, \label{rel3} \\
&\mathcal{C}_n^{(j)}\mathcal{C}_n^{(k)}
=-
\frac{z_j}{z_k}
\mathcal{C}_n^{(k)}\mathcal{C}_n^{(j)}, \ \ \ (j \neq k)
\label{rel4}.
\end{align}
\end{lemma}

\begin{proof}
We show by induction on $n$.  For $n=1$,  from \eqref{initialCD}
$\mathcal{D}_1$ is diagonal and one can directly
see that the relations are satisfied.
For $n$, we assume that $\mathcal{D}_n$ is diagonalizable and write the 
corresponding diagonal matrix as $\mathscr{D}_n=G_n^{-1}\mathcal{D}_n G_n$. 
Also writing  $\mathscr{C}_n=G_n^{-1} \mathcal{C}_n G_n$ and 
$\mathscr{C}_n=\sum_{j=1}^{n} \mathscr{C}_n^{(j)}$, and  noting that
the algebraic relations above do not depend on the choice of basis, we suppose by the
induction hypothesis that the same relations are satisfied by $\mathscr{D}_n$
and $\mathscr{C}_n^{(j)}$. 

We show that the relations hold for $n+1$. To this end, we first
construct $G_{n+1}$. Noting from \eqref{reop1} that $\mathcal{D}_{n+1}$ is an 
upper triangular block matrix whose block diagonal elements are written in 
terms of $\mathcal{D}_n$, 
we assume that $G_{n+1}$ is written as
\begin{equation}
G_{n+1}=
\begin{pmatrix}
G_n &  0 \\
G_n H_n  & G_n
\end{pmatrix},
\label{G-matrix}
\end{equation}
where $2n\times 2n$ matrix $H_n$ remains to be determined. 
Using the induction hypothesis for $n$, one obtains
\begin{align}
&G_{n+1}^{-1}\mathcal{D}_{n+1} G_{n+1} \nonumber \\
=&
\begin{pmatrix}
t \mathscr{D}_n & 0 \\
z_{n+1} \mathscr{D}_n H_n
+(1+t)z_{n+1} \mathscr{C}_n
-t H_n \mathscr{D}_n
& z_{n+1} \mathscr{D}_n
\end{pmatrix}.
\end{align}
The above matrix is guaranteed to be diagonal when 
\begin{equation}
z_{n+1} \mathscr{D}_n H_n
+(1+t)z_{n+1} \mathscr{C}_n
-t H_n \mathscr{D}_n=0.
\end{equation}
Utilizing the above relation and  recalling  $\mathscr{D}_n$
and $\mathscr{C}^{(j)}_n$ satisfy the relation same as that in \eqref{rel2}, 
one finds
\begin{align}
H_n=\mathscr{D}^{-1}_n\sum_{j=1}^n
\frac{(1+t)z_{n+1}}
        {z_j-z_{n+1}} \mathscr{C}_n^{(j)}.
\label{H-matrix}
\end{align}
One thus obtains the diagonal matrix $\mathscr{D}_{n+1}$:
\begin{align}
\mathscr{D}_{n+1}=
\begin{pmatrix}
t \mathscr{D}_n & 0 \\
0 & z_{n+1} \mathscr{D}_n
\end{pmatrix}.
\label{D-matrix}
\end{align}
The remaining task is to derive  $\mathscr{C}_{n+1}^{(j)}$ and
to prove the relations \eqref{rel2}--\eqref{rel4} hold for $n+1$.
Combining  \eqref{reop2}, \eqref{G-matrix} and \eqref{H-matrix},
and also inserting the relations \eqref{rel3} and \eqref{rel4},
one arrives at $\mathscr{C}_{n+1}=\sum_{j=1}^{n+1}\mathscr{C}_{n+1}^{(j)}$
where
\begin{align}
\mathscr{C}_{n+1}^{(j)}=
\begin{cases} \displaystyle
\frac{z_j+t z_{n+1}}{z_j-z_{n+1}}
\begin{pmatrix}
\mathscr{C}_n^{(j)} & 0 \\
0 & -\frac{z_{n+1}}{t} \mathscr{C}_n^{(j)}
\end{pmatrix}  \\[6mm]
\text{ for $1\le j \le n$} \\[6mm]
\begin{pmatrix}
0  & \mathscr{D}_n \\
0 & 0
\end{pmatrix}   \text{ for $j=n+1$}
\end{cases}.
\label{C-matrix}
\end{align}

Finally recalling that $\mathscr{D}_n$ and $\mathscr{C}_n^{(j)}$ 
are supposed to
satisfy the relations \eqref{rel2}--\eqref{rel4} and using the explicit
form of $\mathscr{D}_{n+1}$ \eqref{D-matrix} and $\mathscr{C}_{n+1}^{(j)}$ 
\eqref{C-matrix}, one sees they satisfy the same algebraic relations as those 
in \eqref{rel2}--\eqref{rel4} for $n+1$.
\end{proof}

Due to the algebraic relations \eqref{rel2}, \eqref{rel3}
and \eqref{rel4} in Lemma~\ref{algebra}, 
the matrix product form for the wavefunction \eqref{reov} can be rewritten
as

\begin{align}
&\bra \Phi(\{ z \}_N)| \overline{x_1} \cdots \overline{x_N} \ket
\nonumber \\
=&\prod_{j=1}^N \Bigg( \frac{t}{z_j} \Bigg)^j
\Tr_{W^{\otimes N}}\left[
\mathcal{D}_N^{M-N}
\mathcal{C}_N^{(N)}
\dots\mathcal{C}_N^{(1)} P \right]
 \sum_{\sigma \in S_N} (-1)^\sigma
    \prod_{j=1}^N
  \Bigg( \frac{z_{\sigma(j)}}{t} \Bigg)^{\overline{x_j}} \nonumber \\
=&\prod_{j=1}^N \Bigg( \frac{t}{z_j} \Bigg)^{j-1}
\Tr_{W^{\otimes N}}\left[
\mathcal{D}_N^{M-N}
\mathcal{C}_N^{(N)}
\dots\mathcal{C}_N^{(1)} P \right]
\mathrm{det}_N \Bigg( \Bigg( \frac{z_j}{t} \Bigg)^{\overline{x_k}-1} \Bigg)
\nonumber \\
=&(-1)^{N(N-1)/2} \prod_{j=1}^N \Bigg( \frac{t}{z_j} \Bigg)^{j-1}
\Tr_{W^{\otimes N}}\left[
\mathcal{D}_N^{M-N}
\mathcal{C}_N^{(N)}
\dots\mathcal{C}_N^{(1)} P \right]
\mathrm{det}_N
\Bigg( \Bigg( \frac{z_j}{t} \Bigg)^{\overline{\lambda_{k}}+N-k} \Bigg),
\label{wavefunctiontochu}
\end{align}
where we have used the translation rule
$\overline{\lambda_j}=\overline{x_{N-j+1}}-N+j-1$
between the hole configuration and the Young diagram.
One easily notes that \eqref{wavefunctiontochu}
can be further rewritten in terms of the Schur polynomials:
\begin{align}
\bra \Phi(\{ z \}_N)| \overline{x_1} \cdots \overline{x_N} \ket
=&K s_{\overline{\lambda}} \Bigg( \Bigg\{ \frac{z}{t} \Bigg\}_N \Bigg),
\label{predet}
\end{align}
where the prefactor $K$ given below remains to be determined:
\begin{align}
K=
\prod_{j=1}^N \Bigg( \frac{t}{z_j} \Bigg)^{j-1}
\prod_{1 \le j < k \le N} \frac{z_k-z_j}{t}
\Tr_{W^{\otimes N}}\left[
\mathcal{D}_N^{M-N}
\mathcal{C}_N^{(N)}
\dots\mathcal{C}_N^{(1)} P \right] .
\end{align}

In \eqref{predet},
we notice that the information of the hole configuration
$\{\overline{x_1}, \overline{x_2},\dots,\overline{x_N} \}$
is encoded in the determinant,
while the overall factor $K$ is independent of the configuration.
This fact means that one can determine the factor $K$ by evaluating
the overlap for a particular hole configuration. In fact,
we find the following explicit expression
for the case $\overline{x_j}=j$ ($1\le j \le N$):

\begin{proposition} \label{stepproposition}
The wavefunction $\bra \Phi(\{ z \}_N)| \overline{x_1} \cdots \overline{x_N} \ket$
for the case $\overline{x_j}=j$ ($1\le j \le N$)
has the following form:
\begin{align}
\bra 
\Phi(\{ z \}_N)| \overline{x_1}=1, \cdots, \overline{x_N}=N
\ket
=t^{N(M-N)}
\prod_{1 \le j<k \le N}(z_j+tz_k). \label{stepoverlap}
\end{align}
\end{proposition}

\begin{proof}
We can easily show by its graphical description that
$\bra 
\Phi(\{ z \}_N)| \overline{x_1}=1, \cdots, \overline{x_N}=N
\ket$
can be factorized as
\begin{align}
\bra 
\Phi(\{ z \}_N)| \overline{x_1}=1, \cdots, \overline{x_N}=N
\ket
=t^{N(M-N)}
Z_N(\{ z \}_N), \label{stepoverlapfactorization}
\end{align}
where $Z_N(\{ z \}_N)$ is the domain wall boundary partition function
on an $N \times N$ lattice.
The domain wall boundary partition function
on an $M \times M$ lattice is defined as
\begin{align}
Z_M(\{ z \}_M)=\langle 1 \cdots M|
B(z_1) \cdots B(z_M)| \Omega \rangle,
\end{align}
where $M$ $B$-operators are inserted between the vacuum vector
$| \Omega \rangle=|0 \rangle_1 \otimes \dots \otimes |0 \rangle_M$
and the dual state occupied by particles
$\langle 1 \cdots M|=
{}_1\bra 1|\otimes\dots \otimes{}_M\bra 1|$.

The domain wall boundary partition function
can be analyzed by generalizing it
to the one with inhomogeneties introduced in the quantum spaces
\begin{align}
Z_M(\{ z \}_M|\{ v \}_M)=\langle 1 \cdots M|
B(z_1|\{ v \}_M) \cdots B(z_M|\{ v \}_M)| \Omega \rangle,
\end{align}
where
\begin{align}
B(z|\{ v \}_M)
={}_a \langle 0 |
L_{a M}(z/v_M) \cdots
L_{a 1}(z/v_1)| 1 \rangle_a.
\end{align}
One can show the following factorization formula
for the inhomogeneous domain wall boundary partition function.
\begin{lemma}{\rm cf. \cite{FCWZ}}
\label{inhomogeneousdomainwall}
The domain wall boundary partition function with inhomogeneities
has the following form
\begin{align}
Z_M(\{ z \}_M|\{ v \}_M)
=\prod_{j=1}^M \frac{1}{v_{M-j}}\prod_{1 \le j<k \le M}
(z_j+tz_k). \label{inhomogeneousdomain}
\end{align}
\end{lemma}
Lemma \ref{inhomogeneousdomainwall} can be proved by
using the Izergin-Korepin technique \cite{KBI,Ko,Iz}, i.e.,
show that both hand sides of \eqref{inhomogeneousdomain}
satisfy the same recursive relation, initial condition
and the degree counting of polynomials.

Taking the homogeneous limit $v_j \to 1$ ($j=1,\cdots,M$)
of \eqref{inhomogeneousdomain},
we have \begin{align}
Z_M(\{ z \}_M)
=\prod_{1 \le j<k \le M} (z_j+tz_k). \label{DW}
\end{align}
Replacing $M$ by $N$ in \eqref{DW}
and inserting into \eqref{stepoverlapfactorization},
we have
\begin{align}
\bra 
\Phi(\{ z \}_N)| \overline{x_1}=1, \cdots, \overline{x_N}=N
\ket
=t^{N(M-N)}
\prod_{1 \le j<k \le N}(z_j+tz_k),
\end{align}
which is exactly \eqref{stepoverlap}.
\end{proof}
Using Proposition \ref{stepproposition} into
\eqref{predet}, one can see that
the prefactor $K$ in \eqref{predet}
is determined by the special case of
the dual wavefunction
$\bra 
\Phi(\{ z \}_N)| \overline{x_1}=1, \cdots, \overline{x_N}=N
\ket$
as
\begin{align}
K=\bra 
\Phi(\{ z \}_N)| \overline{x_1}=1, \cdots, \overline{x_N}=N
\ket
=t^{N(M-N)}
\prod_{1 \le j<k \le N}(z_j+tz_k). \label{expressionK}
\end{align}
From \eqref{predet} and \eqref{expressionK},
we have \eqref{dualwavefunctionsandschur},
hence Theorem \ref{theoremdualandschur} is proved.

\end{proof}

\section{Combinatorial formula}
In the previous two sections,
we showed the correspondence between
the dual wavefunction and the Schur polynomials
by giving two different proofs.
To derive a dual Tokuyama-type
combinatorial formula for the Schur polynomials,
one needs to investigate the microscopic structure
and find the partition function expression
for the dual wavefunction.

\begin{figure}[h]
\includegraphics[width=15cm]{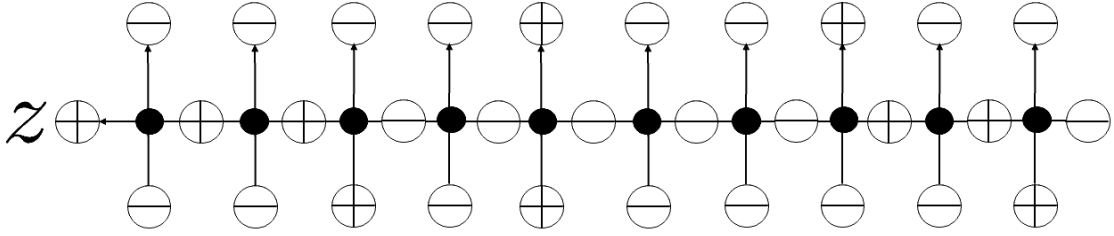}
\caption{The matrix element $\langle \overline{x_1} \cdots \overline{x_N} |B(z)| \overline{y_1}
\cdots \overline{y_{N+1}} \rangle$ for $M=10$, $N=2$,
$\overline{x}=(3,6)$ and
$\overline{y}=(1,6,8)$. One sees that the inner state
is uniquely fixed, and the matrix element is
calculated by multiplying the matrix elements of the $L$-operators
$t \times t \times 1 \times z \times z \times z \times z
\times (t+1)z \times t \times 1=(t+1) t^3 z^5$.}
\label{picturematrixelement}
\end{figure}

The essential thing to find the expression is to view
the dual wavefunction as an object constructed from
$N$ layers of $B$-operators, and analyze the
matrix elements of a single $B$-operator.
One can show the following formula.
\begin{proposition} \label{propositionmatrixelements}
The matrix elements of a single $B$-operator is given by
\begin{align}
\langle \overline{x_1} \cdots \overline{x_N} |B(z)| \overline{y_1}
\cdots \overline{y_{N+1}} \rangle
=&(t+1)^{|\{ \overline{x_j}, \ j=1,\cdots,N \ | \ 
\overline{x_j} \neq \overline{y_{j}}, \ \overline{x_j}
\neq \overline{y_{j+1}} \}|} \nonumber \\
&\times t^{\sum_{j=1}^{N+1} \mathrm{Max}(\overline{x_{j}}
-\overline{y_j}-1, \ 0)}
z^{\sum_{j=1}^N(\overline{y_{j+1}}-\overline{x_j})},
\label{matrixelementscoordinates}
\end{align}
for hole configurations $\{ \overline{x} \}$
$(1 \le \overline{x_1} < \cdots < \overline{x_N} \le M)$ and
$\{ \overline{y} \}$
$(1 \le \overline{y_1} < \cdots < \overline{y_{N+1}} \le M)$
satisfying the interlacing relation
$\overline{y_1} \le \overline{x_1} \le \overline{y_2} \le \overline{x_2} \le \cdots \le \overline{x_N} \le \overline{y_{N+1}}$,
and 0 otherwise.
Here we also set $\overline{x_{N+1}}=M+1$.

Translating into the language of Young diagram via
$\overline{\lambda_j}=\overline{x_{N-j+1}}-N+j-1$, $j=1,\dots,N$,
$\overline{\mu_j}=\overline{x_{N-j+2}}-N+j-2$, $j=1,\dots,N+1$ and
also setting $\overline{\lambda_0}=M-N$, one gets
the following formula for the nonzero matrix elements
when the interlacing relation
 $0 \le \overline{\mu_{N+1}} \le \overline{\lambda_N} \le \overline{\mu_N}
\le \overline{\lambda_{N-1}} \le \cdots \le
\overline{\lambda_1} \le \overline{\mu_1} \le M-N$ is satisfied
\begin{align}
\langle \overline{x_1} \cdots \overline{x_N} |B(z)| \overline{y_1}
\cdots \overline{y_{N+1}} \rangle
=&(t+1)^{|\{ \overline{\lambda_j}, \ j=1,\cdots,N \ | \ 
\overline{\lambda_j} \neq \overline{\mu_{j+1}}, \ \overline{\lambda_j}
\neq \overline{\mu_j}+1 \}|} \nonumber \\
&\times t^{\sum_{j=1}^{N+1} \mathrm{Max}(\overline{\lambda_{j-1}}
-\overline{\mu_j}-1, \ 0)}
z^{\sum_{j=1}^N(\overline{\mu_j}-\overline{\lambda_j})+N}.
\label{matrixelements}
\end{align}

\end{proposition}
\begin{proof}
Let us first count the powers of the spectral parameter $z$.
If the hole configurations
$\{ \overline{x} \}$ and $\{ \overline{y} \}$ are fixed and satisfies
the interlacing relation
$\overline{y_1} \le \overline{x_1} \le \overline{y_2} \le \overline{x_2} \le \cdots \le \overline{x_N} \le \overline{y_{N+1}}$,
the inner states in the auxiliary space is fixed uniquely,
which is a sequence of $0$'s and $1$'s.
We observe that for each sequence $01 \cdots 10$ of the inner states
in the auxiliary space,
all the matrix elements of the
$L$-operators \eqref{loperator} in between contribute to the power $z$,
and gives $z^{\sum_j(\overline{y_{j+1}}-\overline{x_j})}$ for
some sum over $j$. Taking all of the $01 \cdots 10$ sequences into account,
we have the factor $z^{\sum_{j=1}^N(\overline{y_{j+1}}-\overline{x_j})}$.

Let us turn to count the powers of $t+1$ and $t$.
We get a factor $t+1$ for each case when both
$\overline{x_j} \neq \overline{y_j}$ and
$\overline{x_j} \neq \overline{y_{j+1}}$ are satisfied
since the matrix element of the $L$-operator is $[L(z,t)]_{01}^{10}=(t+1)z$
at the $\overline{x_j}$-th site for this case.
One gets $(t+1)^{|\{ \overline{x_j}, \ j=1,\cdots,N \ | \ 
\overline{x_j} \neq \overline{y_{j}}, \ \overline{x_j}
\neq \overline{y_{j+1}} \}|}$ in total.

Next, we count the powers of $t$.
If $\overline{y_j}<\overline{x_j}$ is satisfied,
the matrix elements of the $L$-operators are all
$[L(z,t)]_{01}^{01}=t$ from the $(\overline{y_j}+1)$-th site to 
the $(\overline{x_j}-1)$-th site. On the other hand, $[L(z,t)]_{01}^{01}$
does not appear if $\overline{y_j}=\overline{x_j}$, and there is no
contribution to the power of $t$ for this case.
The contributions from $t$ is given by
$t^{\sum_{j=1}^{N+1} \mathrm{Max}(\overline{x_{j}}
-\overline{y_j}-1, \ 0)}$.

Having calculated all factors, one finds the matrix elements
are given by
\eqref{matrixelementscoordinates} in the coordinate representation.
Translating into the language of Young diagram,
we get \eqref{matrixelements}.
\end{proof}
{\bf Example (coordinate representation)}
Let $M=10$, $N=2$, $\overline{x}=(3,6)$ and
$\overline{y}=(1,6,8)$.
We also set $\overline{x_3}=10+1=11$.
From $\mathrm{Max}(\overline{x_1}-\overline{y_1}-1,0)=\mathrm{Max}(3-1-1,0)=1$, $\mathrm{Max}(\overline{x_2}-\overline{y_2}-1,0)=\mathrm{Max}(6-6-1,0)=0$, $\mathrm{Max}(\overline{x_3}-\overline{y_3}-1,0)=\mathrm{Max}(11-8-1,0)=2$, we have the factor $t^{2+0+1}=t^3$.
The relations $\overline{y_1} \neq \overline{x_1} \neq \overline{y_2}$,
$\overline{y_2} = \overline{x_2} \neq \overline{y_3}$
give the factor $(t+1)^1=t+1$, and we also have the factor $z^5$
from $(\overline{y_2}-\overline{x_1})+
(\overline{y_3}-\overline{x_2})=(6-3)+(8-6)=3+2=5$.
In total, the right hand side of \eqref{matrixelementscoordinates}
is calculated as $(t+1) t^3 z^5$.
One can check that this matches the left hand side
of \eqref{matrixelementscoordinates}, i.e.,
the matrix elements of the corresponding $B$-operator
by explicit calculation
(see Figure \ref{picturematrixelement} for a graphical description
of the corresponding matrix element). \\
\\
{\bf Example (Young diagram representation)}
Let $M=10$, $N=2$, $\overline{x}=(3,6)$ and
$\overline{y}=(1,6,8)$.
We have $\overline{\lambda}=(6-2,3-1)=(4,2)$
and $\overline{\mu}=(8-3,6-2,1-1)=(5,4,0)$.
We also set $\overline{\lambda_0}=10-2=8$.
From $\mathrm{Max}(\overline{\lambda_0}-\overline{\mu_1}-1,0)=\mathrm{Max}(8-5-1,0)=2$, $\mathrm{Max}(\overline{\lambda_1}-\overline{\mu_2}-1,0)=\mathrm{Max}(4-4-1,0)=0$, $\mathrm{Max}(\overline{\lambda_2}-\overline{\mu_3}-1,0)=\mathrm{Max}(2-0-1,0)=1$, we have the factor $t^{2+0+1}=t^3$.
The relations $\overline{\mu_1}+1 \neq \overline{\lambda_1}=\overline{\mu_2}$,
$\overline{\mu_2}+1 \neq \overline{\lambda_2} \neq \overline{\mu_3}$
give the factor $(t+1)^1=t+1$, and we also have the factor $z^5$
from $(\overline{\mu_1}-\overline{\lambda_1})+
(\overline{\mu_2}-\overline{\lambda_2})+2=(5-4)+(4-2)+2=5$.
Altogether, the right hand side of \eqref{matrixelements}
is calculated as $(t+1) t^3 z^5$. \\
\\
In order to descibe the microscopic structure,
we introduce the following strict dual Gelfand-Tsetlin patterns
\begin{eqnarray}
\overline{\mathcal{T}}=\left\{
\begin{array}{ccccccc}
& & & \overline{a_{N-1,N-1}} & & & \\
& & \iddots & & \ddots & & \\
& \overline{a_{1,1}} & & \cdots & & \overline{a_{1,N-1}} & \\
\overline{a_{0,0}} & & \overline{a_{0,1}} & \cdots & \overline{a_{0,N-2}}
& & \overline{a_{0,N-1}}
\end{array}
\right\},
\end{eqnarray}
in which the rows interlace $\overline{a_{i-1,j-1}} \ge \overline{a_{i,j}}
\ge \overline{a_{i-1,j}}$,
and the entries in horizontal rows are strictly decreasing.
We use this picture so as to be consistent with the dual wavefunction
description.
The strict Gelfand-Tsetlin patterns essentially label
the inner states by recording the positions
of the holes in the auxiliary spaces.

For each strict Gelfand-Tsetlin pattern, we assign
the following weight:
\begin{align}
\overline{G}(\overline{\mathcal{T}},\{ z \}_N)
=\prod_{i=0}^{N-1}
\prod_{j=i}^{N-1} \overline{\gamma}(\overline{a_{i,j}})
z_N^{d_0(\overline{\mathcal{T}})-d_1(\overline{\mathcal{T}})}
z_{N-1}^{d_1(\overline{\mathcal{T}})-d_2(\overline{\mathcal{T}})} \cdots
z_2^{d_{N-2}(\overline{\mathcal{T}})-d_{N-1}(\overline{\mathcal{T}})}
z_1^{d_{N-1}(\overline{\mathcal{T}})}, \label{productofweights}
\end{align}
where
$d_j(\overline{\mathcal{T}})=\sum_{k=j}^{N-1} \overline{a_{j,k}}$, $j=0,\dots,N-1$
is the sum of the entries of the Gelfand-Tsetlin pattern
in the $j$-th row counted from the bottom,
and $\overline{\gamma}(\overline{a_{i,j}})$ is defined as
\begin{align}
\overline{\gamma}(\overline{a_{i,j}})
=t^{\mathrm{Max}(\overline{a_{i+1,j}}-\overline{a_{i,j}}-1,0)} \times
\left\{
\begin{array}{ll}
t+1 & \overline{a_{i-1,j}} \neq \overline{a_{i,j}} \neq \overline{a_{i-1,j-1}}
\\
1 & \mathrm{otherwise}
\end{array}
\right.
, \label{bargamma}
\end{align}
for pairs of integers $(i,j)$ satisfying
$0 \le i \le N-1, \ i \le j \le {N-1}$.
Note that we define $\overline{\gamma}(\overline{a_{0,j}})$, $j=0, \dots,N-1$
since we need these weights to describe the dual wavefunction
and the dual Tokuyama-type formula
(whereas one does not need to define $\gamma(a_{0,j})$ to describe
the original wavefunction and the Tokuyama formula).
We also define $\overline{a_{j,j-1}}=M$ for $j=1,\cdots,N$.

As again, the inner states making non-zero contributions
can be characterized by the strict Gelfand-Tsetlin pattern
with the bottom row fixed by the Young diagram as
$\overline{a_{0,j}}=\overline{\lambda_{j+1}}+N-j-1$.

\begin{theorem} \label{anothercombinatorialformula}
We have the following combinatorial formula for the Schur polynomials
\begin{align}
&t^{N(M-N)} \prod_{1 \le j < k \le N}(z_j+t z_k)
s_{\overline{\lambda}} \Bigg( \Bigg\{ \frac{z}{t} \Bigg\}_N \Bigg) \nonumber \\
=&\sum_{\overline{\mathcal{T}}}
\overline{G}(\overline{\mathcal{T}},\{ z \}_N) \nonumber \\
=&\sum_{\overline{\mathcal{T}}}
\prod_{i=0}^{N-1}
\prod_{j=i}^{N-1} \overline{\gamma}(\overline{a_{i,j}})
z_N^{d_0(\overline{\mathcal{T}})-d_1(\overline{\mathcal{T}})}
z_{N-1}^{d_1(\overline{\mathcal{T}})-d_2(\overline{\mathcal{T}})} \cdots
z_2^{d_{N-2}(\overline{\mathcal{T}})-d_{N-1}(\overline{\mathcal{T}})}
z_1^{d_{N-1}(\overline{\mathcal{T}})}, \label{newcombinatorialformula}
\end{align}
where the sum is over all strict Gelfand-Tsetlin patterns
with the bottom row of the strict Gelfand-Tsetlin pattern
is fixed by the Young diagram as
$\overline{a_{0,j}}=\overline{\lambda_{j+1}}+N-j-1$.
\end{theorem}

\begin{proof}
The proof follows from evaluating the dual wavefunction \\
$\langle \Phi(z_1,\dots,z_N)| \overline{x_1} \cdots \overline{x_N}
\rangle=
\langle 1 \cdots M|B(z_1) \cdots B(z_N)| \overline{x_1} \cdots \overline{x_N}
\rangle$ in two ways. First, we note from Theorem \ref{theoremdualandschur}
that $\langle 1 \cdots M|B(z_1) \cdots B(z_N)| \overline{x_1} \cdots \overline{x_N} \rangle$ is expressed using Schur polynomials as
\eqref{dualwavefunctionsandschur}.

Another way of evaluation can be accomplished by viewing the dual wavefunction
as a partition function constructed from $N$ layers of $B$-operators,
inserting the completeness relation
and decomposing it as
sums of products of matrix elements of the $B$-operators.
That is, we decompose the dual wavefunction as
\begin{align}
&\langle 1 \cdots M|B(z_1) \cdots B(z_N)| \overline{x_1} \cdots \overline{x_N}
\rangle \nonumber \\
=&\sum_{\overline{\mathcal{T}}}
\langle 1 \cdots M | \nonumber \\
&\times \prod_{j=1}^N
\Bigg\{ |\overline{a_{N-j+1,N-1}}+1, \cdots, \overline{a_{N-j+1,N-j+1}}+1
\rangle
\langle \overline{a_{N-j+1,N-1}}+1, \cdots, \overline{a_{N-j+1,N-j+1}}+1 |
B(z_j) \Bigg\} \nonumber \\
&\times | \overline{a_{0,N-1}}+1 \cdots \overline{a_{0,0}}+1 \rangle,
\label{intermediatewavefunction}
\end{align}
where
\begin{align}
| \overline{a_{j,N-1}}+1 \cdots \overline{a_{j,j}}+1 \rangle
=\prod_{k=j}^{N-1}
\sigma_{\overline{a_{j,k}}+1}^+(|1 \rangle_1 \otimes \cdots \otimes |1 \rangle_M)
, \\
\langle \overline{a_{j,N-1}}+1 \cdots \overline{a_{j,j}}+1|=
(\langle 1|_1 \otimes \cdots \otimes \langle 1|_M) \prod_{k=j}^{N-1}
\sigma_{\overline{a_{j,k}}+1}^-.
\end{align}
and use the formula for the matrix elements of a single $B$-operator
\eqref{matrixelementscoordinates}
in Proposition \ref{propositionmatrixelements}.
Then one finds the product of the matrix elements of the $B$-operators
in \eqref{intermediatewavefunction}
corresponding to each strict Gelfand-Tsetlin pattern $\overline{\mathcal{T}}$
can be expressed as
\eqref{productofweights}
\begin{align}
\overline{G}(\overline{\mathcal{T}},\{ z \}_N)
=\prod_{i=0}^{N-1}
\prod_{j=i}^{N-1} \overline{\gamma}(\overline{a_{i,j}})
z_N^{d_0(\overline{\mathcal{T}})-d_1(\overline{\mathcal{T}})}
z_{N-1}^{d_1(\overline{\mathcal{T}})-d_2(\overline{\mathcal{T}})} \cdots
z_2^{d_{N-2}(\overline{\mathcal{T}})-d_{N-1}(\overline{\mathcal{T}})}
z_1^{d_{N-1}(\overline{\mathcal{T}})}.
\end{align}
Hence, the identity
\eqref{intermediatewavefunction}
can be rewritten in the following form
for the dual wavefunction
\begin{align}
&\langle 1 \cdots M|B(z_1) \cdots B(z_N)| \overline{x_1} \cdots \overline{x_N}
\rangle \nonumber \\
=&\sum_{\overline{\mathcal{T}}}
\overline{G}(\overline{\mathcal{T}},\{ z \}_N) \nonumber \\
=&\sum_{\overline{\mathcal{T}}}
\prod_{i=0}^{N-1}
\prod_{j=i}^{N-1} \overline{\gamma}(\overline{a_{i,j}})
z_N^{d_0(\overline{\mathcal{T}})-d_1(\overline{\mathcal{T}})}
z_{N-1}^{d_1(\overline{\mathcal{T}})-d_2(\overline{\mathcal{T}})} \cdots
z_2^{d_{N-2}(\overline{\mathcal{T}})-d_{N-1}(\overline{\mathcal{T}})}
z_1^{d_{N-1}(\overline{\mathcal{T}})}. \label{forcomparison}
\end{align}
Comparing the two expressions
\eqref{dualwavefunctionsandschur}
and \eqref{forcomparison} evaluated by two ways,
we get \eqref{newcombinatorialformula}.

\end{proof}

The combinatorial formula
\eqref{newcombinatorialformula}
in Theorem \ref{anothercombinatorialformula}
can be rewritten into the following form by scaling every
spectral parameter $z_j$ to $tz_j$
and cancelling powers of $t$ of both hand sides and simplyfing.
\begin{theorem} \label{dualcombinatorialformula}
We have the following combinatorial formula for the Schur polynomials
\begin{align}
&\prod_{1 \le j < k \le N}(z_j+t z_k)
s_{\overline{\lambda}}(\{ z \}_N) \nonumber \\
=&t^{\sum_{j=1}^N \overline{\lambda_j}-N(M-N)}\sum_{\overline{\mathcal{T}}}
\prod_{i=0}^{N-1}
\prod_{j=i}^{N-1} \overline{\gamma}(\overline{a_{i,j}})
z_N^{d_0(\overline{\mathcal{T}})-d_1(\overline{\mathcal{T}})}
z_{N-1}^{d_1(\overline{\mathcal{T}})-d_2(\overline{\mathcal{T}})} \cdots
z_2^{d_{N-2}(\overline{\mathcal{T}})-d_{N-1}(\overline{\mathcal{T}})}
z_1^{d_{N-1}(\overline{\mathcal{T}})}, \label{newcombinatorialformulaanotherexpression}
\end{align}
where the sum is over all strict Gelfand-Tsetlin patterns
with the bottom row of the strict Gelfand-Tsetlin pattern
is fixed by the Young diagram as
$\overline{a_{0,j}}=\overline{\lambda_{j+1}}+N-j-1$.
\end{theorem}

Let us discuss the differences between
Theorem \ref{combinatorialformula} and Theorem \ref{dualcombinatorialformula}.
In the original wavefunction, the factor
$\gamma(a_{i,j})$ \eqref{gamma}
in \eqref{combinatorialformulaequation}
depends only on three neighbors $a_{i,j}$, $a_{i-1,j}$ and $a_{i-1,j-1}$
in the strict Gelfand-Tsetlin pattern.
On the other hand, for the case of the dual wavefunction,
the factor
$\overline{\gamma}(\overline{a_{i,j}})$ \eqref{bargamma} in
\eqref{newcombinatorialformulaanotherexpression} depends on
four neighbors $\overline{a_{i,j}}$, $\overline{a_{i-1,j}}$, $\overline{a_{i-1,j-1}}$ and $\overline{a_{i+1,j}}$.
Note also the order of the symmetric variables (spectral parameters)
in \eqref{combinatorialformulaequation} is $z_1,\dots,z_N$,
while it is $z_N,\dots,z_1$ in
\eqref{newcombinatorialformulaanotherexpression}, i.e., the order is reversed.
Moreover,
the right hand side of \eqref{newcombinatorialformulaanotherexpression}
has powers of $t$ as factors,
which explicitly depends on the size of the Young diagram,
the total number of sites $M$ and the total number of holes $N$.
This explicit dependence cannot be found in
\eqref{combinatorialformulaequation}. \\
\\
{\bf Example}
Let us check \eqref{newcombinatorialformulaanotherexpression}
by an example.
Consider the case
$M=4$, $N=2$, $\overline{\lambda}
=(\overline{\lambda}_1,\overline{\lambda_2})=(2,1)$.
The bottom row of the dual strict Gelfand-Tsetlin pattern
is fixed as $\overline{a_{0,0}}=\overline{\lambda_1}+N-0-1=2+2-0-1=3$,
$\overline{a_{0,1}}=\overline{\lambda_2}+N-1-1=1+2-1-1=1$.
From the interlacing relation $3=\overline{a_{0,0}} \ge \overline{a_{1,1}} \ge \overline{a_{0,1}}=1$, we have $\overline{a_{1,1}}=1$, $2$ or $3$.
Therefore, there are three strict Gelfand-Tsetlin
patterns in the sum of
\eqref{newcombinatorialformulaanotherexpression}

\begin{eqnarray}
&\overline{\mathcal{T}}=\left\{
\begin{array}{ccc}
 & \overline{a_{1,1}} & \\
\overline{a_{0,0}} &  & \overline{a_{0,1}}
\end{array}
\right\} \nonumber \\
=&\left\{
\begin{array}{ccc}
 & 1 & \\
3 &  & 1
\end{array}
\right\}, \
\left\{
\begin{array}{ccc}
 & 2 & \\
3 &  & 1
\end{array}
\right\}, \
\left\{
\begin{array}{ccc}
 & 3 & \\
3 &  & 1
\end{array}
\right\}.
\end{eqnarray}
We also set $\overline{a_{1,0}}=\overline{a_{2,1}}=4$
for each strict Gelfand-Tsetlin pattern.
Keeping this in mind,
let us calculate the weights for each pattern. \\
\\
(1) $\displaystyle
\overline{\mathcal{T}}=\left\{
\begin{array}{ccc}
 & 1 & \\
3 &  & 1
\end{array}
\right\}.
$ \\
We have
$\overline{\gamma}(\overline{a_{0,0}})=t^{\mathrm{Max}(4-3-1,0)}=1$,
$\overline{\gamma}(\overline{a_{0,1}})=t^{\mathrm{Max}(1-1-1,0)}=1$,
$\overline{\gamma}(\overline{a_{1,1}})=t^{\mathrm{Max}(4-1-1,0)}=t^2$,
$d_0(\overline{\mathcal{T}})=3+1=4$, $d_1(\overline{\mathcal{T}})=1$.
Thus we have $\overline{\gamma}(\overline{a_{0,0}})
\overline{\gamma}(\overline{a_{0,1}}) \overline{\gamma}(\overline{a_{1,1}})
z_2^{d_0(\overline{\mathcal{T}})-d_1(\overline{\mathcal{T}})}
z_1^{d_1(\overline{\mathcal{T}})}=t^2 z_2^3 z_1$.
\\
\\
(2) $\displaystyle
\overline{\mathcal{T}}=\left\{
\begin{array}{ccc}
 & 2 & \\
3 &  & 1
\end{array}
\right\}.
$ \\
We have
$\overline{\gamma}(\overline{a_{0,0}})=t^{\mathrm{Max}(4-3-1,0)}=1$,
$\overline{\gamma}(\overline{a_{0,1}})=t^{\mathrm{Max}(2-1-1,0)}=1$,
$\overline{\gamma}(\overline{a_{1,1}})=t^{\mathrm{Max}(4-2-1,0)}(t+1)
=t(t+1)$,
$d_0(\overline{\mathcal{T}})=3+1=4$, $d_1(\overline{\mathcal{T}})=2$.
The corresponding weight is \\
$\overline{\gamma}(\overline{a_{0,0}})
\overline{\gamma}(\overline{a_{0,1}}) \overline{\gamma}(\overline{a_{1,1}})
z_2^{d_0(\overline{\mathcal{T}})-d_1(\overline{\mathcal{T}})}
z_1^{d_1(\overline{\mathcal{T}})}=t(t+1)z_2^2 z_1^2$.
\\
\\
(3) $\displaystyle
\overline{\mathcal{T}}=\left\{
\begin{array}{ccc}
 & 3 & \\
3 &  & 1
\end{array}
\right\}.
$ \\
We have
$\overline{\gamma}(\overline{a_{0,0}})=t^{\mathrm{Max}(4-3-1,0)}=1$,
$\overline{\gamma}(\overline{a_{0,1}})=t^{\mathrm{Max}(3-1-1,0)}=t$,
$\overline{\gamma}(\overline{a_{1,1}})=t^{\mathrm{Max}(4-3-1,0)}=1$,
$d_0(\overline{\mathcal{T}})=3+1=4$, $d_1(\overline{\mathcal{T}})=3$.
We have
$\overline{\gamma}(\overline{a_{0,0}})
\overline{\gamma}(\overline{a_{0,1}}) \overline{\gamma}(\overline{a_{1,1}})
z_2^{d_0(\overline{\mathcal{T}})-d_1(\overline{\mathcal{T}})}
z_1^{d_1(\overline{\mathcal{T}})}=tz_2 z_1^3$ in total. \\
\\
Summing the three weights calculated above and noting
$\sum_{j=1}^2 \overline{\lambda_j}-N(M-N)=2+1-4=-1$,
the right hand side of
\eqref{newcombinatorialformulaanotherexpression}
is
$t^{-1}(t^2 z_2^3 z_1+t(t+1)z_2^2 z_1^2+tz_2 z_1^3)
=tz_2^3 z_1+(t+1)z_2^2 z_1^2+z_2 z_1^3$,
which can be factorized as $(z_1+tz_2)(z_1^2 z_2+z_1 z_2^2)
=(z_1+tz_2)s_{(2,1)}(z_1,z_2)$,
which is exactly the left hand side of
\eqref{newcombinatorialformulaanotherexpression}.

\section{A generalization of the correspondence}
We have showed Theorem \ref{theoremdualandschur}
which gives the relation between the dual wavefunction
and the Schur polynomials, for which we gave two proofs.
The one given in section 4 can be applied to a generalization
of the Felderhof model, where inhomogeneous parameters
are now introduced in the quantum spaces.
Since the original wavefunction was found to give the
factorial Schur polynomials \cite{BMN}, one expects
the dual wavefunction also gives the factorial Schur polynomials.

The $L$-operator which constructs the wavefunction
now has dependence on the quantum space $\mathcal{F}_j$: at the $j$-th site
in the quantum space, we introduce the following $L$-operator
\begin{eqnarray}
L_{aj}(z,t,\alpha_j)=\left( 
\begin{array}{cccc}
1 & 0 & 0 & 0 \\
0 & t & 1 & 0 \\
0 & (t+1)z & z+\alpha_j & 0 \\
0 & 0 & 0 & z-t \alpha_j
\end{array}
\right). \label{generalizedfelderhofloperator}
\end{eqnarray}
The $L$-operators $L_{aj}(z,t,\alpha_j)$
now has inhomogeneous parameters
$\alpha_j$, $j=1,\cdots M$ besides the spectral parameter
and the deformation parameters.

We consider the wavefunction
$\langle x_1 \cdots x_N|\Psi(z_1,\dots,z_N,\{ \alpha \}) \rangle$
by introducing the $N$-particle state
\begin{align}
\Psi(z_1,\dots,z_N,\{ \alpha \}) \rangle
=
B(z_1,\{ \alpha \}) \cdots
B(z_N,\{ \alpha \})| \Omega \rangle,
\end{align}
where the $B$-operator
\begin{align}
B(z,\{ \alpha \})
={}_a \langle 0 |
L_{a M}(z,t,\alpha_M) \cdots
L_{a 1}(z,t,\alpha_1)| 1 \rangle_a,
\end{align}
now has dependence on the inhomogeneous parameters
$\{ \alpha \}=\{ \alpha_1,\dots,\alpha_M \}$,
which turns out to be the factorial parameters
of the factorial Schur polynomials defined below.

\begin{definition}
The factorial Schur polynomials is defined
to be the following determinant:
\begin{align}
s_\lambda(\{ z \}_N|\{ \alpha \}|)
=\frac{F_{\lambda+\delta}(\{ z \}_N|\{ \alpha \})}
{\prod_{1 \le j < k \le N}(z_j-z_k)}, \label{factorialSchurfunction}
\end{align}
where $\{ z \}=\{z_1,\dots,z_N \}$ is a set of variables
and $\lambda$ denotes a Young diagram
$\lambda=(\lambda_1,\lambda_2,\dots,\lambda_N)$
with weakly decreasing non-negative integers
$\lambda_1 \ge \lambda_2 \ge \cdots \ge \lambda_N \ge 0$,
and $\delta=(N-1,N-2,\dots,0)$.
$F_{\mu}(\{ z \}_N|\{ \alpha \})$
is an $N \times N$ determinant
\begin{align}
F_{\mu}(\{ z \}_N|\{ \alpha \})
=\mathrm{det}_N
\Bigg(
\prod_{j=1}^\mu(z_k+\alpha_j)
\Bigg).
\end{align}
We remark that one must respect the ordering of the factorial parameters
$\{ \alpha \}=\{\alpha_1,\dots,\alpha_M \}$.
\end{definition}
Bump, McNamara and Nakasuji showed the following correspondence
between the wavefunction of the Felderhof model with
inhomogeneties and the factorial Schur polynomials.
\begin{theorem} \cite{BMN}
The wavefunction
$\langle x_1 \dots x_N|\Psi(z_1,\dots,z_N,\{ \alpha \}) \rangle$
is expressed by the factorial Schur polynomials as
\begin{align}
\langle x_1 \dots x_N|\Psi(z_1,\dots,z_N,\{ \alpha \}) \rangle
=\prod_{1 \le j<k \le N}(z_j+tz_k)
s_\lambda(\{ z \}_N|\{ \alpha \}),
\end{align}
under the relation $\lambda_j=x_{N-j+1}-N+j-1$, $j=1,\dots,N$.
\end{theorem}
This Theorem was proved by noting that the arguments in \cite{BBF}
naturally lift to this inhomogeneous setting.
One first shows that the wavefunction is a polynomial
of $t$ with highest weight degree $N(N-1)/2$.
Then one evaluates the wavefunction at $t=-1$,
at which the six-vertex model reduces to a five-vertex model,
and each configuration making nonzero contribution
to the wavefunction essentially corrresponds to each term
of the determinant expansion of the numerator of the
factorial Schur polynomials \eqref{factorialSchurfunction}.

Let us now state the result for the dual wavefunction
$\langle \Phi(z_1,\dots,z_N,\{ \alpha \})|
\overline{x_1} \cdots \overline{x_N} \rangle$
which is the overlap between the hole configurations
$|\overline{x_1} \cdots \overline{x_N} \rangle$
and the dual $N$-particle state
\begin{align}
\langle \Phi(z_1,\dots,z_N,\{ \alpha \})|
=\langle 1 \cdots M|
B(z_1,\{ \alpha \}) \cdots
B(z_N,\{ \alpha \}).
\end{align}

By applying the argument in section 4, one gets the
following relation between the dual wavefunction and the
factorial Schur polynomials.
\begin{theorem} \label{theoremdualandfactorialschur}
The dual wavefunction
$\langle \Phi(z_1,\dots,z_N,\{ \alpha \})| \overline{x_1} \cdots \overline{x_N}
\rangle$
can be expressed by the factorial Schur polynomials as
\begin{align}
\langle \Phi(z_1,\dots,z_N,\{ \alpha \})| \overline{x_1} \cdots \overline{x_N}
\rangle=t^{N(M-N)} \prod_{1 \le j < k \le N}(z_j+t z_k)
s_{\overline{\lambda}} \Bigg( \Bigg\{ \frac{z}{t} \Bigg\}_N,
\{-\alpha \} \Bigg).
\label{dualwavefunctionsandfactorialschur}
\end{align}
Here the Young diagram for the factorial Schur polynomials
correspond to the particle configuration under the relation
$\overline{\lambda_j}=\overline{x_{N-j+1}}-N+j-1$, $j=1,\dots,N$,
and the symmetric variables are
$\displaystyle \Bigg\{ \frac{z}{t} \Bigg\}_N=
\Bigg\{ \frac{z_1}{t},\dots,\frac{z_N}{t} \Bigg\}$.
Moreover, the signs of the parameters of the factorial Schur polynomials
in the right hand side of \eqref{dualwavefunctionsandfactorialschur} are now
inverted simultaneously: $\{-\alpha \}=\{-\alpha_1,\dots,-\alpha_M \}$.
\end{theorem}
The correspondence \eqref{dualwavefunctionsandfactorialschur}
includes the special case $t=1$
of the relation between the dual wavefunction and factorial Schur polynomials
in \cite{BMN},
which was proved by starting from the result for the relation
between the original wavefunction and the factorial Schur polynomials,
using arguments on the symmetry of the $L$-operators
to transform the original correspondence to the dual correspondence.
This argument seems very difficult for the case $t \neq 1$
even for the ordinary Schur polynomials.
However, one can naturally lift the arguments given in section 4
to this inhomogeneous setting. The problem reduces to 
the case of the $t=-1$, where the six-vertex model reduces to
the five-vertex model. Since we now have the inhomogenous parameters,
this introduction of additional parameters is reflected in the final
expression of the correspondence in
\eqref{dualwavefunctionsandfactorialschur}.

\section{Conclusion}
We investigated the Felderhof free-fermion model,
and analyzed the dual wavefunction in two ways.
We first showed the precise relation between the dual
wavefunction and the Schur polynomials,
in which we gave two proofs in sections 4 and 5 respectively.
One by using the arguments
by \cite{BBF}, and another one by combining the
matrix product method and the analysis on
the domain wall boundary partition function.
Next, by calculating the matrix elements
of a single $B$-operator, we give a combinatorial
expression of the Schur polynomials in terms of
strict Gelfand-Tsetlin patterns.
By comparing the two expressions,
we obtained a combinatorial formula
of the Schur polynomials, which can be regarded as
a dual version of the Tokuyama formula,
since it was found \cite{BBF} that the original wavefunction naturally gives
a realization of the Tokuyama formula for the Schur polynomials,
and we are now dealing with the dual wavefunction.

We also generalized the relation between the dual wavefunction
to the Felderhof model with inhomogeneous parameters
in the quantum space and the factorial Schur polynomials,
which is motivated by the fact that the
wavefunction of the Felderhof model with inhomogeneties
are given by the factorial Schur polynomials \cite{BMN}.
The expression can be extended furthermore
to the Felderhof model with two types of
inhomogeneous parameters, and there are
correspondences between the original and the dual
wavefunctions and a generalization
of the factorial Schur polynomials \cite{MSW1,MSW2}.

One of the important problems related to
this paper is to study
the dual wavefunction for the case of other boundary conditions
and find
combinatorial formulas for other symmetric polynomials
such as the symplectic Schur and Schur $Q$ functions.
See \cite{Iv,BBCG,Tabony,HK,BS} for examples
for the relation with the
wavefunctions and the Felderhof model with other boundary conditions.

The Schur polynomials appears not only as the wavefunction
of the Felderhof model, but also as special limits
of the wavefunction XXZ-type six-vertex model.
The integrable five-vertex model which is the $t=0$ limit
of the $L$-operator \eqref{loperator},
which gives the Schur polynomials, can be regarded as
special limits of both the Felderhof model and the XXZ model.
See \cite{MS2,MSW1,BW,BWZ,Korff,GK,GK2} for examples on the
recent investigations on the
combinatorics of the symmetric polynomials
from the viewpoint of partition functions,
in which the combinatorial identities of
various symmetric polynomials such as the 
Schur, Grothendieck, Hall-Littlewood and their
noncommutative versions are derived.

We finally remark that in recent works, it is revealed by number theorists
that the six-vertex model considered in this paper
can be regarded as a special case of the ``metaplectic ice",
which is a six-vertex model over a non-archimedean local field
(see \cite{Meta} for example).
It seems worthwhile to study these models
and find novel combinatorial formulas
by means of modern statistical physical methods
and techniques developed to analyze quantum integrable models.

\section*{Acknowledgements}
Thus work was partially supported by grant-in-Aid
for Research Activity start-up No. 15H06218
and Scientific Research (C) No. 16K05468.

\end{document}